\documentclass{elsarticle}

\usepackage{lineno,hyperref}
\usepackage[]{graphicx}
\usepackage[]{color}
\usepackage{amsmath}
\usepackage{amsfonts}
\usepackage{mathrsfs}
\usepackage{amssymb}
\usepackage{amsbsy}
\usepackage{multicol,multirow}
\usepackage{bigints}
\usepackage{algorithm}
\usepackage{algpseudocode}
\usepackage{rotating}
\usepackage{url}
\modulolinenumbers[5]

\definecolor{pinegreen}{rgb}{0,0.55,0.3}

\definecolor{orange}{rgb}{1,0.5,0}

\definecolor{dkgreen}{rgb}{0,0.5,0}

\newcommand{\bR}{\mathbb{R}}
\newcommand{\bZ}{\mathbb{Z}}
\newcommand{\expe}{\ensuremath{\mathrm{e}}} 
\newcommand{\supp}{\ensuremath{\mathrm{supp}}} 
\newcommand{\ie}{\textit{i.e.}}
\newcommand{\eg}{\textit{e.g.}}
\newcommand {\mcF}{\mathcal{F}} 
\newcommand {\mcW}{\mathcal{W}} 
\newcommand {\bS}{\mathbb{S}}   
\newcommand {\msSH}{\mathscr{S\!\!H}} 
\newcommand {\mcR}{\mathcal{R}}  
\newcommand {\mcT}{\mathcal{T}}  
\newcommand {\mcP}{\mathcal{P}}  

\ifpdf 
	\DeclareGraphicsExtensions{{.pdf}}
	\graphicspath{{../NETNAmanuscript_APNUMsub_gen2016/IMACS_RecImages/PDF/}}
\else 
	\DeclareGraphicsExtensions{{.eps}}
	\graphicspath{{../NETNAmanuscript_APNUMsub_gen2016/IMACS_RecImages/EPS/}}
\fi

\journal{Applied Numerical Mathematics}

\begin{document}

\begin{frontmatter}

\title{Numerical assessment of shearlet-based regularization in ROI 
tomography}

\author[FerraraAddress]{Tatiana A.~Bubba\corref{CorrespAuthor}}
\cortext[CorrespAuthor]{Corresponding author}
\ead{bbbtnl@unife.it}

\author[HoustonAddress]{Demetrio Labate}
\ead{dlabate@math.uh.edu}

\author[FerraraAddress]{Gaetano Zanghirati}
\ead{gaetano.zanghirati@unife.it}

\author[FerraraAddress]{Silvia Bonettini}
\ead{silvia.bonettini@unife.it}

\address[FerraraAddress]{Dept. of Mathematics and Computer Science, University of Ferrara, and INdAM-GNCS, via G.~Saragat 1, 44122 Ferrara, Italy}
\address[HoustonAddress]{Department of Mathematics, University of Houston, 651 Phillip G.~Hoffman, Houston, TX USA 77204-3008}

\begin{abstract}
When it comes to computed tomography (CT), the possibility to reconstruct a small region-of-interest (ROI) using truncated projection data is particularly appealing due to its potential to lower radiation exposure and reduce the scanning time. However, ROI reconstruction from truncated projections is an ill-posed inverse problem, with the ill-posedness becoming more severe when the ROI size is getting smaller. To address this problem, both \textit{ad hoc} analytic formulas and iterative numerical schemes have been proposed in the literature. In this paper, we introduce a novel approach for ROI CT reconstruction, formulated as a convex optimization problem with a regularized term based on shearlets. Our numerical implementation consists of an iterative scheme based on the scaled gradient projection (SGP) method and is tested in the context of fan beam CT. Our results show that this approach is essentially insensitive to the location of the ROI and remains very stable also when the ROI size is rather small.
\end{abstract}


\end{frontmatter}

\linenumbers

\section{Introduction}
\label{sec:intro}  

Computed Tomography (CT) is a noninvasive imaging technique designed 
to visualize the internal structure of a body or an object without surgical 
intervention or destructive material testing. To generate CT images, 
X-rays are propagated through the object and \emph{projections} are 
collected from multiple views so that the density of the object can be 
reconstructed solving an appropriate inverse problem. 
The impact of CT has been enormous in areas including industrial 
nondestructive testing, security screening and medical diagnostics even 
though, in this last application, exposure to X-ray radiation comes with 
health hazards for patients. Luckily, in many biomedical situations, such as 
contrast-enhanced cardiac imaging or some surgical implant procedures, 
one is interested in examining only a small region-of-interest (ROI) with high 
resolution. Hence, there is no need to irradiate the entire body but only a 
smaller region, with the advantage of reducing radiation exposure and 
shortening the scanning time.

However, CT reconstruction is an ill-posed problem and the ill-posedness is 
even more severe when one attempts to solve the reconstruction problem 
from truncated projections \cite{natterer:tomography}, as in the case of ROI 
CT.  Indeed, the direct application of classical reconstruction algorithms such 
as the Filtered Back-Projection (FBP) or the FDK algorithms 
\cite{natterer:imagerec} (with the missing projection data set to zero) 
typically produces unacceptable visual artifacts, and is more and more 
unstable to noise as the size of the ROI decreases.

During the last decade, a number of analytic and algebraic methods has 
been proposed to address the problem of ROI reconstruction from truncated 
projections \cite{clackdoyle:roi}.  In particular, it was shown that it is possible 
to derive analytic ROI reconstruction formulas from truncated data, even 
though such formulas usually require restrictive assumptions on the location 
of the ROI and depend on the acquisition setting 
\cite{noo:imrec,clackdoyle:alcif, zou:image}. For example, the Differentiated 
Back-Projection (DBP) methods \cite{noo04,Clackdoyle04,zeng:iterative} 
can be applied only if there exists at least one projection view in which the 
complete (\ie{}, non-truncated) projection data are available. 
Algebraic or iterative methods, on the other hand, are generally more
flexible, since they can be applied to essentially any type of acquisition 
mode, also including several physical processes in the modeling. However, 
they are usually computationally more demanding, even if advances in 
high-performance computing make algebraic methods more and more 
competitive \cite{Beister12,Guorui:katsevichGPU}.
Some common iterative algorithms, that have been adapted to ROI 
tomography, are the simultaneous iterative reconstruction technique (SIRT) 
\cite{Herman76}, the maximum likelihood expectation-maximization 
algorithm (MLEM) \cite{Shepp82} and the least-squares Conjugate Gradient 
(LSCG) method \cite{Hestenes52}.
They typically involve some form of prior knowledge on the object 
attenuation function, or a regularization term to ensure a stable ROI 
reconstruction \cite{Hamelin_2010,YangYuJiang_2010,Ziegler2008}. 
Anyway, the performance of existing methods is usually rather sensitive to 
the ROI size and the presence of noise.

In this paper, we address the ROI CT reconstruction problem using an
iterative minimization algorithm formulated as a convex optimization problem with a regularized functional based on 
shearlets, a multiscale methods especially designed to sparsely approximate images with 
edges~\cite{ELL08}. The approach here we consider relies on part on the setup recently
proposed in~\cite{Goossens14}, where the reconstruction method combines a data fidelity constraint
to ensure that reconstructed projections matched the observed data inside the ROI and a data consistency condition to enforce that solutions are consistent with the tomographic reconstruction problem. Starting from this setup, we propose two formulations of the regularized reconstruction problem that we call implicit and explicit formulations. One main novelty of this paper   
is that our numerical implementation is based on the scaled gradient projection 
method~\cite{Bonettini09}, which is an accelerated first-order descent 
method for convex and non-convex objective functions. This method is particularly 
effective when the feasible region is given by ``simple'' constraints,  that is, 
when projecting onto such a feasible region is not a heavy task. 
The numerical tests reported in this paper show that our algorithm produces accurate ROI 
reconstructions for any ROI location and also for ROI sizes that are small 
with respect to the field of view. We remark that a preliminary version of these results were recently presented in a conference paper~\cite{Bubba15}.

The paper is organized as follows: in Section \ref{sec:continuous} we 
introduce the problem formulation in both the continuous and the discrete 
settings, translating it into two possible formulations of the 
objective function. Also, we quickly reiterate the mathematical background 
for both the shearlets and the SGP algorithm.
In Section \ref{sec:distancedriven} we describe the setup for the simulated 
data and we recall the basics on distance-driven method \cite{DeMan04}.
Numerical experiments and results are given in Section \ref{sec:NumEx}. 
Finally, we draw some conclusions in Section \ref{sec:conclusions}.

\section{SGP-shearlet approach for the region-of-interest tomography problem} 
\label{sec:continuous} 

In this section, we introduce the ROI CT problem in two dimensions. Our 
framework applies to different projection geometries, including parallel and 
fan-beam, and extends to higher dimensions. However, for brevity, we only 
consider here the two-dimensional case. We derive also a regularized 
reconstruction algorithm using two different formulations for the objective 
function. Our implementation is an iterative approach based on a first order 
descent method that can be applied to reach the unique solution of the 
problem.

\subsection{Continuous setup: Radon transform and ROI solution uniqueness}
\label{sec:contsetup}

The CT reconstruction problem consists in reconstructing a density function from a set of projections, 
obtained by measuring attenuation over straight lines.
Mathematically, it becomes
a ``line-integral model" through the Radon transform notion.
Given a function $f \in L^1(\bR^2)$, the \emph{Radon transform} 
of $f$ at $(\theta,\tau)$ is the line integral of $f$ over the lines (or rays) 
$\ell(\theta,\tau)$ perpendicular to $\mathbf{e}_{\theta} = \bigl(\cos(\theta),\sin(\theta)\bigr)^T \in \bS^1$ with (signed) distance $\tau \in \bR$ from the origin:
\[
\mcR f(\theta,\tau) = \int_{\ell (\theta,\tau)} f(\mathbf{x}) \, d\mathbf{x} = 
\int_{\bR^2} \delta(\tau - \mathbf{x} \cdot \mathbf{\expe}_{\theta}) \, f(\mathbf{x}) \, d\mathbf{x},
\]
where $\ell(\theta,\tau) = \{ \mathbf{x} \in \bR^2 : \mathbf{x} \cdot \mathbf{e}_\theta = \tau \}$.
Therefore, the Radon transform maps $f$ into the set of its linear projections defined on the tangent space
$\mcT = \{ (\theta,\tau) : \theta \in [0,2\pi), \, \tau \in \bR \}$.
In what follows, we will
refer to the Radon projections as the \emph{full sinogram} and we shall
denote it by:
\begin{equation}
y(\theta, \tau) = \mcR f(\theta,\tau), \quad \theta \in [0,2\pi), \, \tau \in \bR.
\label{eq:RadonProj}
\end{equation}
We will refer to the domain of \eqref{eq:RadonProj} as the 
\emph{projection domain} and to the domain of the density function as the 
\emph{image domain}.

In the ROI tomographic problem, \textit{projections are only collected for those rays meeting a region of interest} inside the field of view. The goal is to recover the density function inside the ROI, while the rest of the function is essentially ignored.
Denoting the ROI as $S \subset \bR^2$, the set of ROI-truncated projections is identified to be the set
\[
\mcP(S) = \{ (\theta, \tau) \in \mcT : \ell(\theta,\tau) \cap S \neq  \emptyset \} \subset \mcT.
\]
Thus, the ROI reconstruction problem can be formulated as the problem of reconstructing the function $f$  restricted to the ROI $S$ from the truncated Radon projections:
\begin{equation}\label{eq:ROIeq}
y_0(\theta, \tau) = M(\theta, \tau) \, \mcR f(\theta, \tau),
\qquad
M(\theta, \tau) = 1_{\mcP(S)}(\theta, \tau)
\end{equation}
where $M$ is the \emph{mask function} 
and $1_A$ is the indicator function of the set $A$. We will refer to $y_0$ 
as the \emph{truncated sinogram}.
In the following, we will assume that the ROI is a disk $S \subset \bR^2$ with 
center $\mathbf{p}_{\text{ROI}} \in \bR^2$ and radius $R_{\text{ROI}} \in 
\bR$. Thus, $\mcP(S) = |\tau - \mathbf{p}_{\text{ROI}} \cdot 
\mathbf{e}_\theta| < R_{\text{ROI}}$. Clearly, more general convex ROIs 
can be handled by finding the minimal enclosing disk for this ROI and 
reconstructing the image for this disk.

A natural approach for obtaining a stable reconstruction of $f$ from equation 
\eqref{eq:ROIeq} is by computing the least squares solution $\widehat{f}$
\begin{equation}
\widehat{f} = \arg \min_{f} \| M \mcR f - y_0 \|_2^2.
\end{equation}
However, the minimizer of this problem is not unique in general, since the
set of solutions is the affine subspace
\[
V = \{ f \in L^2(\bR^2) : y = \mcR f \quad \text{and} \quad My=y_0 \}.
\] 
Indeed, it is known that the solution of the ROI problem, in general, is not 
guaranteed to be unique \cite{natterer:tomography}. Even when the 
uniqueness is ensured, the inversion of the Radon transform is an ill-posed 
problem and the ill-posedness may be more severe when projections are 
truncated, as in the case of the ROI CT problem. 
A classical approach to achieve uniqueness is by using Tikhonov 
regularization, which imposes an additional norm condition by searching for 
the minimum-norm solution. This norm condition can be applied in 
the image domain or in the projection domain, leading to different solutions, in general.

To define a norm condition in the projection domain, we can proceed as in~\cite{Goossens14}
where the Riesz operator is used to relate the norm of $f$ to the norm of $\mcR f$. Namely,
recall the definition of the Riesz potential operator 
$I^{-\alpha}$ for a function $g \in L^2(\mcT)$:
\begin{equation}\label{eq:Riesz}
\mcF(I^{\alpha} g) (\theta,\xi) = |\xi|^{-\alpha} \mcF g(\theta,\xi), 
\qquad \alpha < 2.
\end{equation}
Recall also the following formula \cite{natterer:imagerec}:
\begin{equation}\label{eq:RieszEquiv}
f = \frac{1}{4\pi} \mcR^* I^{-1} \mcR f,
\end{equation}
where $\mcR^*$ is the adjoint operator of $\mcR$, also known as 
backprojection operator. Hence, we have:
\begin{equation}\label{eq:RieszNorm}
\begin{split}
\| f \|_2^2 \; = \; \frac{1}{4\pi} \langle f, \mcR^* I^{-1} \mcR f \rangle  
\; &= \;  \frac{1}{4\pi} \langle \mcR f, I^{-1} \mcR f \rangle \\
 &= \;  \frac{1}{4\pi} \langle I^{-\frac{1}{2}} \mcR f, I^{-\frac{1}{2}} \mcR f \rangle 
\; = \;  \frac{1}{4\pi} \left\| I^{-\frac{1}{2}} \mcR f \right\|_2^2.
\end{split}
\end{equation}
Thus, we are lead to the following optimization problem:
\begin{equation}\label{eq:MinProbl}
\widehat{f} = \arg \min_{f} 	\; \bigl\|I^{-1/2} \mcR f \bigr\|_2^2
\qquad \quad \text{s. t.} 
\quad y = \mcR f \quad \text{and} \quad My = y_0.
\end{equation}

To formulate an objective function for the solution of the ROI reconstruction problem, we adopt the following observations from \cite{Goossens14}. 
It is easy to see that
\begin{equation}\label{eq:inpainting}
y = y_0 + (1-M) \, y,
\end{equation}
indicating that the ROI reconstruction problem can be viewed 
as an extrapolation problem where $y_0$, given on $\mcP(S)$, is to be extrapolated 
outside $\mcP(S)$. Clearly we cannot choose \emph{any} extrapolated
function $y$ outside $\mcP(S)$ but the following global constraint is needed
\begin{equation}\label{eq:GlobCon}
Rf = y.
\end{equation}
Equation \eqref{eq:GlobCon} ensures that $y$ belongs to range of the Radon transform of the density  
function $f \in L^1(\bR) \cap L^2(\bR)$. By applying $M$ and $1 - M$ to the 
left-hand 
and right-hand sides of (\ref{eq:GlobCon}), respectively, we hence obtain the 
following equations (cf.~\cite{Goossens14}): 
\begin{align}\label{eq:Fid}
M \mcR f &= M \, y = y_0	\qquad\; \text{(data fidelity)} \\
(1-M) \mcR f  &= (1 - M) \, y \qquad \text{(data consistency)} \label{eq:Const}
\end{align}
The data fidelity equation defines a constraint inside the ROI and the data 
consistency equation enforces accurate reconstruction inside the ROI.
These equations can be 
combined with regularization to obtain a computationally suitable algorithm 
for the ROI CT problem, as we shall see in section \ref{sec:DiscreteObjFun}.

\subsection{Discrete framework: explicit and implicit formulation of the objective function} 
\label{sec:DiscreteObjFun}

To derive a discrete formulation for the ROI CT reconstruction framework, 
we need first to discretize equations \eqref{eq:Fid} and \eqref{eq:Const}. 
We shall denote by $K$ the number of 
projection  angles and by $P$ the number of detector cells (\ie{}, samples 
along the detector array); $N$ is both the width and the height in pixels of the 
object to reconstruct. 
Given a projection geometry, the matrix $\mathbf{W}$ of the forward problem, 
that represents the map from the image domain to the projection domain, 
has dimensions $KP \times N^2$. 
The mask identifying the ROI is the diagonal matrix $\mathbf{M}$ of 
dimensions $KP \times KP$ whose entries are either 0 or 1.
The unknown discrete density function $f$ to be reconstructed is represented as a vector $\mathbf{f}$ of length $N^2$ in which the entries are stacked column by column. 
Similarly, the full sinogram $y$ and the truncated sinogram 
$y_0$ are represented as vectors $\mathbf{y}$ and $\mathbf{y}_0$, 
respectively, of length $KP$, obtained by 
stacking the entries column by column. We recall that $\mathbf{y}$ and 
$\mathbf{y}_0$ are related to each other by the data fidelity equation
\[
\mathbf{y}_0 = \mathbf{M y} = \mathbf{M W f},
\]
and the data consistency equation, that sets the extrapolation scheme 
outside the ROI, reads as
\[
(\mathbf{I}_{KP} - \mathbf{M}) \mathbf{W f} = (\mathbf{I}_{KP} - \mathbf{M})\mathbf{y},
\]
where $\mathbf{I}_{KP}$ is the $KP \times KP$ identity matrix. 
As indicated above, data fidelity and data consistency equations need 
to be coupled with regularization to yield a unique solution. Thus, similarly to 
\eqref{eq:MinProbl}, we obtain the following discrete optimization problem:
\begin{equation}\label{eq:OptProblDiscr}
\widehat{\mathbf{f}} = \arg \min_{\mathbf{f}} \| \mathbf{\Phi W f} \|_2^2
\qquad \quad \mbox{s. t.} 
\quad \mathbf{y} = \mathbf{W f} \quad \mbox{and} \quad 
\mathbf{M y} = \mathbf{y}_0
\end{equation}
where $\mathbf{\Phi}$ is a discrete filter corresponding to the Riesz potential 
operator $I^{-\frac{1}{2}}$.

However, rather than requiring the exact equalities stated by data fidelity 
and data consistency equations, we shall minimize the $L^2$-norm error 
associated to them, according to a maximum likelihood approach.   
By exploiting this idea, we can state two different optimization problems and 
objective functions.
On the one hand, we can consider the $L^2$-norm error of the data fidelity 
equation only and incorporate the data consistency information in the regularization term. In this case, the only
variable we are minimizing on is the image $\mathbf{f}$ to be reconstructed.
This approach yields:
\begin{equation*}
\widehat{\mathbf{f}} = 
\arg \min_{\mathbf{f} \in \Omega_{\mathbf{f}}} \mathbf{\Psi}(\mathbf{f})
\end{equation*}
where
\begin{equation}\label{eq:Implicit}
\mathbf{\Psi}(\mathbf{f}) = \frac{1}{2} \| \mathbf{M W f} - \mathbf{y}_0 \|_2^2 +
\lambda \, \| \mathbf{\Phi} ((\mathbf{I}_{KP} - \mathbf{M}) \mathbf{W f} + 
\mathbf{y}_0) \|_2^2.
\end{equation}
We will refer to \eqref{eq:Implicit} as the \emph{implicit formulation}, since the 
full sinogram $\mathbf{y}$ does not appear explicitly.

On the other hand, we can consider both the data fidelity and the data 
consistency $L^2$-norm errors into the objective function to assess if, by 
explicitly incorporating the extrapolation scheme in the minimization problem,
we achieve a more accurate reconstruction. In this case, the image $\mathbf{f}$ 
to be reconstructed and the full sinogram $\mathbf{y}$ are both unknowns.
This second approach yields:
\[
\Bigl(\widehat{\mathbf{f}}, \widehat{\mathbf{y}}\Bigr) = \arg \min_{\substack{\mathbf{f} \in \Omega_{\mathbf{f}} \\ \mathbf{y} \geq \mathbf{0}}} 
\mathbf{\Psi}(\mathbf{f}, \mathbf{y})
\]
where
\begin{equation}\label{eq:Explicit}
\begin{split}
\mathbf{\Psi}(\mathbf{f}, \mathbf{y}) &= 
\frac{1}{2} \| \mathbf{M W f} - \mathbf{y}_0 \|_2^2 +
\frac{1}{2} \| (\mathbf{I}_{KP} - \mathbf{M}) (\mathbf{W f} - \mathbf{y}) \|_2^2 
\\ 
&\qquad+ \lambda \, \| \mathbf{\Phi} ((\mathbf{I}_{KP} - \mathbf{M}) \mathbf{y} + 
\mathbf{y}_0) \|_2^2.
\end{split}
\end{equation}
In contrast with the above formulation, we will refer to \eqref{eq:Explicit} 
as the \emph{explicit formulation}. Notice that each term of the objective function is convex 
with respect to the unknowns.
In both cases, $\lambda$ denotes the regularization parameter and the 
feasible region $\Omega_\mathbf{f}$ is either defined as $\mathbf{f} \geq 0$ or, when the maximum pixel value $L$ of the image happens to be known, as $0 \leq \mathbf{f} \leq L$. 
Here the inequalities are meant component-wise. Notice also that each term of the objective function is convex with respect to the unknowns.

We will also include a smoothed total variation (TV) term~\cite{Vogel02} in the objective 
function to control oscillatory artifacts in the numerical solution.
In conclusion, our objective functions are
\begin{equation}\label{eq:ImplicitTV}
\mathbf{\Psi}(\mathbf{f}) = \dfrac{1}{2} \| \mathbf{M W f} - \mathbf{y}_0 \|_2^2 +
\lambda \, \| \mathbf{\Phi} ((\mathbf{I}_{KP} - \mathbf{M}) \mathbf{W f} + 
\mathbf{y}_0) \|_2^2 +
\rho \; \text{TV}_\delta (\mathbf{f})
\end{equation}
for the implicit formulation and 
\begin{equation}\label{eq:ExplicitTV}
\begin{split}
\mathbf{\Psi}(\mathbf{f}, \mathbf{y}) &= 
\frac{1}{2} \| \mathbf{M W f} - \mathbf{y}_0 \|_2^2 +
\frac{1}{2} \| (\mathbf{I}_{KP} - \mathbf{M}) (\mathbf{W f} - \mathbf{y}) \|_2^2 
\\ 
&\qquad
+\lambda \, \| \mathbf{\Phi} ((\mathbf{I}_{KP} - \mathbf{M}) \mathbf{y} + 
\mathbf{y}_0) \|_2^2  +
\rho \; \text{TV}_\delta (\mathbf{f})
\end{split}
\end{equation}
for the explicit one.
Here, $\rho$ is a regularization parameter and $\delta$ is the TV smoothing 
parameter.
TV was chosen since it is widely used in medical imaging. 
In the following, we will investigate if using the classical Tikhonov-like 
regularization term is superior to the TV minimization. 

We remark that, in the absence of noise, the solution 
$\widehat{\mathbf{f}}$ is the same for both the 
implicit and the explicit cases (up to a vector in the nullspace of 
$\mathbf{W}$). However, in the presence of noise, the explicit formulation 
should benefit from the presence of the norm 
$\| \mathbf{W} \mathbf{f} -\mathbf{y}\|_2^2$ that is in general non-zero (and can be controlled using an appropriate stopping criterion, typically depending on the noise variance), forcing a more accurate reconstruction \emph{inside} the ROI.

\subsection{Regularization term: discrete shearlet transform}
\label{sec:reguterms}

In the section above, the term $\mathbf{\Phi}$ was intentionally left 
unspecified. Indeed, in place of considering a straightforward 
matrix discretization of the Riesz potential operator, we will 
approximate it by using the discrete shearlet transform, 
a multiscale method which refines the conventional wavelet framework by 
combining multiscale analysis and directional 
sensitivity \cite{GuoLab:2007,KL12_book}. For completeness, we briefly recall below the main ideas about shearlets.

Let $\phi \in L^2(\bR^2)$. A 2D affine family generated by $\phi$ is a 
collection of functions of the form:
\[
\biggl\{ \phi_{M, t} (x) = |\det(M) \, |^{-\frac{1}{2}} \phi \Bigl(M^{-1} (x-t) \Bigr) \; : \; M \in G, t \in \bR^2 \biggr\}
\]
where $G$ is a subset of the group $GL_2(\bR)$ of invertible $2 \times 2$ 
matrices.
The function $\phi$ is called a 2D \emph{continuous wavelet} if 
\[
g(x) = \int_G \int_{\bR^2} \langle g, \, \phi_{M, t} \rangle \; \phi_{M, t}(x) \, dt \, 
d\lambda(M)
\qquad \forall \; g \in L^2(\bR^2),
\]
where $\lambda(M)$ is an suitable measure on $G$. The corresponding 2D 
\emph{continuous wavelet transform} of $g$ is the mapping
\[
g \quad \mapsto \quad \mcW g(M,t) = \langle g, \, \phi_{M, t} \rangle, 
\qquad \quad
M \in G, t \in \bR^2.
\]
Discrete wavelet transforms are obtained by discretizing 
$\mcW g(M,t)$ on an appropriate set. Usually, $M \in G$ and $t \in \bR^2$ 
are replaced by $A^j$, $j \in \bZ$, and $k \in \bZ^2$, respectively, for appropriate choices of the matrix $A$. Shearlets are obtained from the 2D affine systems
\[
\biggl\{ \phi_{a,s,t}(x) = {|\det(M_{a,s}) \, |}^{-\frac{1}{2}} \phi \Bigl(M_{a,s}^{-1} (x-t) \Bigr) \; : \; M \in \Gamma, t \in \bR^2 \biggr\}, 
\]
where $\phi \in L^2(\bR^2)$ and $\Gamma = \left\{ M_{a,s} = \left( \begin{array}{cc}
a & s \sqrt{a} \\ 0 & \sqrt{a}
\end{array} \right) 
\; : \; a \in \bR^+, s \in \bR  \right\}$ is a subset of $GL_2(\bR)$.
Observe that the matrix $M_{a,s}$ is obtained by multiplying the \emph{anisotropic dilation} matrix $A_a$ with the \emph{shear} matrix $S_s$:
\[
M_{a,s} = \left( \begin{array}{cc}
a & s \sqrt{a} \\ 0 & \sqrt{a}
\end{array} \right) = 
\left( \begin{array}{cc}
a & 0 \\ 0 & \sqrt{a}
\end{array} \right)
\left( \begin{array}{cc}
1 & s \\ 0 & 1
\end{array} \right) =: A_a S_s.
\]
Here the variables $a \in \bR^+$, $s \in \bR$ and $t \in \bR^2$ denote the scale, 
orientation and the spatial location, respectively. Thus, shearlets are formed 
by dilating, shearing and translating an appropriate  mother shearlet function 
$\phi \in L^2 (\bR^2)$ \cite{Guo09}. Roughly speaking, they are well localized 
waveforms whose orientation is controlled by the shear parameter $s$ and 
that become increasingly elongated at fine scales (as $a \rightarrow 0$).
The \emph{continuous shearlet transform} of $g$ is the map
\[
g \quad \mapsto \quad \msSH g(a,s,t) = \langle g, \, \phi_{a,s,t} \rangle, 
\qquad \quad
a \in \bR^+, s \in \bR, t \in \bR^2.
\]

Discrete shearlet systems are formally defined by sampling a continuous shearlet systems on an appropriate discrete set. In particular, we choose
\begin{equation}
\begin{gathered}
\biggl\{ \phi_{j,k,m} (x) = 2^{\frac{3}{4}j} \phi \Bigl(S_k A_{2^j}x - m \Bigr)\; : \; j, k \in \bZ, m \in \bZ^2 \biggr\} \\
\text{with $\phi$ s. t.} \quad
\mcF (\phi) (\omega) = \mcF (\phi_1) (\omega_1) \; \mcF (\phi_2) \left( 
\frac{\omega_2}{\omega_1} \right),
\end{gathered}
\label{eq:discrsh}
\end{equation}
where $\omega = (\omega_1, \omega_2)$, $ \mcF (\phi_1)$ is 
the Fourier transform of a wavelet function with compact support away from the origin and $\mcF (\phi_2)$ is a 
compactly supported bump function with $\supp(\mcF (\phi_2)) \subset 
[-1,1]$. 
Such a shearlet is called \emph{classical shearlet}.
Starting from a classical shearlet and under other suitable assumptions (see 
\cite{Guo06}), one can obtain a tight frame for $L^2 (\bR^2)$. 
This indicates that the decomposition is invertible and the transformation is 
numerically well-conditioned. Notice that there exist several choices for 
$\phi_1$ and $\phi_2$ satisfying the classical shearlet definition. One possible choice is to set $\phi_1$ to be a Lemari\`{e}--Meyer wavelet and $\phi_2$ to be a spline bump function.

Summarizing, shearlets exhibit very appealing mathematical properties: they are well localized, namely they are compactly supported in the frequency domain and have fast decay in the spatial domain; the parabolic scaling allows to take care of the anisotropic structures; highly directional sensitivity is provided by the shearing parameter.
Thanks to these properties, 
shearlets provide optimally sparse approximations with images containing 
$C^2$-edges, outperforming conventional wavelets \cite{GuoLab:2007}. 
This is potentially relevant in CT-like applications, since point-like structures 
in the image domain map onto sine-shaped curvilinear structures in the 
projection domain.

\subsection{Scaled gradient projection method}
\label{sec:SGP}

In this section, we introduce our algorithm for the solution 
of the optimization problems \eqref{eq:ImplicitTV}-\eqref{eq:ExplicitTV}, 
whose step-by-step description is given in Algorithm~\ref{alg:SGP}. 
The proposed method, called \textit{scaled gradient projection} (SGP) 
method \cite{Bonettini09}, is an 
iterative approach from the family of first-order descent methods that apply to 
convex (and non-convex), differentiable and constrained problems with 
``simple'' feasible regions. This is indeed the case of our formulation, given 
that the feasible region $\Omega_{\mathbf{f}}$ is usually either a box or a 
non-negativity constraint and, as already indicated, our functionals are 
convex and differentiable.

\begin{algorithm}[t]
\caption{Scaled Gradient Projection Method}\label{alg:SGP}
\begin{algorithmic}
\State Choose the starting point $\mathbf{f}^{(0)} \in \Omega_{\mathbf{f}}$, 
set the parameters $\beta, \gamma \in (0,1)$ and 
$0 < \alpha_{\min} < \alpha_{\max}$. Fix a positive integer $\mu$.
\For {$k = 0, 1, 2, \dots$} 
\State Step 1. Choose the parameter $\alpha_k \in [\alpha_{\min}, \alpha_{\max}]$ and the scaling matrix $D_k \in \mathcal{D}_L$ ;
\State Step 2. Projection: $\mathbf{z}^{(k)} = \mathcal{P}_{\Omega_f} (\mathbf{f}^{(k)} - \alpha_k D_k \nabla \Psi (\mathbf{f}^{(k)}))$. 
\quad  \textbf{if} $\mathbf{z}^{(k)} = \mathbf{f}^{(k)}$ 
           \textbf{then} stop:  
          $\mathbf{f}^{(k)}$ is a stationary point; \textbf{end if}
\State Step 3. Descent direction: $\mathbf{d}^{(k)} = \mathbf{z}^{(k)} - \mathbf{f}^{(k)}$;
\State Step 4. Set $\lambda_k = 1$ and $\Psi_{\max} = \max_{0 \leq j \leq \min(k,\mu-1)} \Psi(\mathbf{f}^{(k-j)})$;
\State Step 5. Backtracking loop:
\If {$\Psi(\mathbf{f}^{(k)} + \lambda_k \mathbf{d}^{(k)}) \leq \Psi_{\max} + \beta \lambda_k \nabla \Psi(\mathbf{f}^{(k)})^T \mathbf{d}^{(k)}$}
\State go to Step 6;	
\Else
\State set $\lambda_k =\gamma \lambda_k$ and go to Step 5;
 \EndIf 
\State Step 6. Set $\mathbf{f}^{(k+1)} = \mathbf{f}^{(k)} + \lambda_k \mathbf{d}^{(k)}$.
\EndFor
\end{algorithmic}
\end{algorithm}

Some remarks about SGP are in order.
When the objective function reads as in \eqref{eq:ImplicitTV}, the $(k+1)$-th 
iteration is 
\begin{equation}
\mathbf{f}^{(k+1)} = (1-\lambda_k) \, \mathbf{f}^{(k)} + 
\lambda_k \, \mathcal{P}_{\Omega_\mathbf{f}}  
\Bigl( \mathbf{f}^{(k)} - \alpha_k D_k \nabla \Psi (\mathbf{f}^{(k)}) \Bigr)
\qquad
k = 0, 1, 2, \ldots,
\end{equation}
where $\lambda_k$, $\alpha_k$ are suitable steplengths, $D_k$ is the 
scaling matrix and $\mathcal{P}_{\Omega_\mathbf{f}}$ is the projector onto 
the feasible region.
The main feature of this method consists in the combination of non-expensive 
diagonally scaled gradient directions with steplength selection rules specially 
designed for these directions.  
In details, any choice of the steplength $\alpha_k$ in a closed interval 
$[\alpha_{\min}, \alpha_{\max}] \subset \bR^+$ and of the scaling matrix 
$D_k$ in the compact set 
$\mathcal{D}_L$ is allowed, where $\mathcal{D}_L$ is the set of the symmetric positive definite matrices $D$ such that 
$\|D\| \leq L$ and $\|D^{-1}\| \leq L$, for a given threshold $L > 1$. 
This is very important from a practical point of view because it makes their 
updating rules problem-related and performance-aware. 
In particular, SGP is equipped with an adaptive steplength selection based on 
the Barzilai-Borwein (BB) updating rules \cite{Barzilai88,Fletcher01}. 
In practice, by means of a variable threshold, one of the two different 
selection strategies
\begin{align*}
\alpha_k^{\text{BB1}} &= \arg \min_{\alpha_k \in \bR} \bigl\| B(\alpha_k) \, 
\mathbf{s}^{(k-1)} - \mathbf{\zeta}^{(k-1)} \bigr\| \\
\alpha_k^{\text{BB2}} &= \arg \min_{\alpha_k \in \bR} \bigl\| \mathbf{s}^{(k-1)} - B(\alpha_k)^{-1} \mathbf{\zeta}^{(k-1)} \bigr\| 
\end{align*}
is selected \cite{Frassoldati08}, 
where $B(\alpha_k) = (\alpha_k D_k)^{-1}$ approximates
the Hessian matrix $\nabla^2 \Psi(\mathbf{f}^{(k)})$,
$\mathbf{s}^{(k-1)} = \mathbf{f}^{(k)} - \mathbf{f}^{(k-1)}$ and 
$\mathbf{\zeta}^{(k-1)} = \nabla \Psi(\mathbf{f}^{(k)}) - \nabla 
\Psi(\mathbf{f}^{(k-1)})$.
As far as the scaling matrix concerns, the updating rule for each entry 
$d_i^{(k)}$ is:
\begin{equation}
d_i^{(k)} = \min \left\{ \sigma, \max \left\{ \frac{1}{\sigma}, \mathbf{f}_i^{(k)} 
\right\} \right\}
\qquad i = 1, \, \ldots, KP
\end{equation}
where $\sigma$ is an appropriate threshold.
We can setup the SGP method for the case of two unknowns, as in equation 
\eqref{eq:ExplicitTV}, by applying the SGP iteration to the ``enlarge'' variable 
$(\mathbf{f}, \mathbf{y})$. This is accomplished by tacking a block-diagonal scaling matrix: the block corresponding to the image $\mathbf{f}$ reads as described above, while for the block corresponding to the full sinogram
$\mathbf{y}$ we take $D_k = \mathbf{I}_{KP}$ for each $k = 0, 1, 2, \ldots$. That is, the $(k+1)$-th iteration, with $k = 0, 1, 2, \ldots$, is
\begin{align}
\mathbf{f}^{(k+1)} &= (1-\lambda_k) \, \mathbf{f}^{(k)} + 
\lambda_k \, \mathcal{P}_{\Omega_\mathbf{f}}  
\Bigl( \mathbf{f}^{(k)} - \alpha_k D_k \nabla_{\mathbf{f}} 
\Psi (\mathbf{f}^{(k)}, \mathbf{y}^{(k)}) \Bigr) \\
\mathbf{y}^{(k+1)} &= (1-\lambda_k) \, \mathbf{y}^{(k)} + 
\lambda_k \, \mathcal{P}_{\mathbf{y} \geq 0}  
\Bigl( \mathbf{y}^{(k)} - \alpha_k \nabla_{\mathbf{y}} 
\Psi (\mathbf{f}^{(k)}, \mathbf{y}^{(k)}) \Bigr)
\end{align}

We notice also that global convergence properties are 
ensured by exploiting a nonmonotone line-search strategy along the feasible 
direction \cite{Grippo86,Dai05}. Indeed, the nonmonotone line-search 
strategy, as implemented in step 5, ensures that $\Psi(\mathbf{f}^{(k+1)})$ 
($\Psi(\mathbf{f}^{(k+1)}, \mathbf{y}^{(k+1)})$, resp.) is lower than the 
maximum of the objective function on the last $\mu$ iterations; of course, if 
$\mu = 1$ then the strategy reduces to the standard monotone Armijo rule.
Moreover, since classical shearlets are tight
frame, the equivalence $\mathbf{\Phi}^T \mathbf{\Phi} = \mathbf{I}_{KP}$ 
holds true and this allows to reduce the computational complexity of each 
algorithm iteration due to simpler updating rules.

Finally, notice that the SGP algorithm can be also used as an iterative 
regularization method applied to the un-regularized functional, by means of 
an early stopping technique.

\section{Data simulation}
\label{sec:distancedriven}

To demonstrate and validate our approach, we use a synthetic data set 
known as \emph{modified Shepp-Logan phantom}. It is available, for instance, in the Matlab Image Processing toolbox.
All phantom data are simulated using the geometry of a micro-CT scanner 
used for real measurements. We simulate 2D fan-beam data over $182$ 
uniformly spaced angles over $2\pi$. The detector consists of $130$ 
elements with a pixel pitch of $0.8 mm$. The distance between tube and 
detector is set to $291.20 mm$ and the radius of rotation is $115.84 mm$. 
The detector is offset by $1.5$ pixels. The 2D fan-beam geometry, 
represented by the matrix $\mathbf{W}$ introduced in section 
\ref{sec:DiscreteObjFun}, is implemented by using the state-of-the-art 
method called distance-driven \cite{DeMan04}. A brief introduction to its main 
mathematical ideas is reported in the section \ref{sec:dd}.

\subsection{Distance-driven method}
\label{sec:dd}

We recall that $\mathbf{W}$ is the $KP \times N^2$ matrix mapping the image domain into the projection domain. 
Every row of $\mathbf{W}$ contains weights that relate the pixel intensities in 
the image domain to the corresponding sample in the projection domain. 
Obviously, the value of each weight depends on the chosen interpolation 
scheme. 
The distance-driven method combines a highly 
sequential memory access pattern with relatively low arithmetic complexity, 
without introducing interpolation artifacts in the image or projection domains.
Essentially, this approach is based on converting the projection problem into 
a 1D re-sampling problem. 
There are two main ingredients in the distance-driven method. The first one is 
the \emph{kernel operation}:
\begin{equation}\label{eq:DistanceDriven}
\begin{gathered}
b_n = \sum_j  w_j c_j 
\qquad \text{with} \\
w_j = \dfrac{\bigl[ \min(\xi_{m+1}, \upsilon_{n+1}) - \max (\xi_m, \upsilon_n) \bigr]_+}{\upsilon_{n+1} - \upsilon_n} \,, \quad
[x]_+ = \max(x,0) \,,
\end{gathered}
\end{equation}
which allows one to compute the destination signal values $\{b_j\}_j$ from the 
sample values $\{ c_i \}_i$ of a source signal, the sample source locations 
$\{ \xi_i\}_i$ and the sample destination locations $\{ \upsilon_j \}_j$.
The second element of the method is that there is  a (possibly zero) length of overlap between each image pixel and each detector cell due to the bijection between the position on the 
detector and the position within an image row (or column). Thus, every point 
within an image row is uniquely mapped onto a point on the detector, and vice versa. 
In practice, to compute the overlap length, all pixel boundaries in an image 
row and all detector cell boundaries are mapped onto a common line, \eg{}, 
a line parallel to a coordinate axis. That is, each length of overlap is the 
interval length between two adjacent intersections, obtained by connecting 
the boundaries midpoints of all detector cells and pixels in a image row to 
the X-ray source and by computing the intercepts of these lines with the 
common axis.
The final weights are achieved by normalizing the overlap length by the detector cell width. This corresponds 
exactly to apply the kernel operation from equation \eqref{eq:DistanceDriven}.
In our case, $b_n$ is the theoretical (unblurred and noiseless) value 
measured at the $n$-th detector cell and $c_j$ is the estimate of the $j$-th 
pixel attenuation function. For example, if the $n$-th detector cell is ``shadowed'' by only two pixels in a row, equation \eqref{eq:DistanceDriven} 
reads as:
\[
b_n = \frac{\xi_{m+1} - \upsilon_n}{\upsilon_{n+1} - \upsilon_n} c_m+ \frac{\upsilon_{n+1} - \xi_{m+1}}{\upsilon_{n+1} - \upsilon_n} c_{m+1}.
\]
Notice that an efficient implementation of the distance-driven method can 
be achieved by using a vector-oriented approach instead of processing one 
pixel at a time. Indeed, one can update an entire row of the matrix (\ie{}, 
chosen one detector cell at a fixed angular view) in one shot by identifying 
\emph{all} the pixels that shadow that detector cell and their corresponding 
length of overlap. Moreover, the detector can be chosen as the common 
line hence avoiding the need to calculate the projections of the detector cell 
boundaries. Finally, when the common line is not the detector, an effective choice is to 
adaptively select the common line by mapping onto the $x$-axis or $y$-axis, depending on the angular view.

\section{Numerical experiments}
\label{sec:NumEx}

\begin{table}[t]
\centering
\makebox[0pt][c]{%
\begin{tabular}{@{}c|r|rrl@{\ }|@{\ }rrl@{\ }|@{\ }rrl@{}}
\multicolumn{2}{c|}{}& \multicolumn{3}{c}{$\gamma = 0.5 N$} & \multicolumn{3}{c}{$\gamma = 0.3 N$} & \multicolumn{3}{c}{$\gamma = 0.15 N$}  \\
\multicolumn{1}{c}{}&& value & iter & param & value & iter & param & value & iter & param  \\
\hline
\multirow{5}{*}{\rotatebox{90}{PSNR}}&
SGP + Sh + TV  & 51.35 & 1277 & $\lambda=5\cdot 10^{-4}$ 
& 40.35 & 1951 & $\lambda=5\cdot 10^{-4}$
& 41.04 & 432 & $\lambda=5\cdot 10^{-4}$ \\
&& & & $\rho=0.01$ & & & $\rho=1$ & & & $\rho=1$  \\
&SGP + Sh & 40.66 & 426 & $\lambda=5\cdot 10^{-4}$ 
& 35.85 & 745 & $\lambda=5\cdot 10^{-4}$
& 29.56 & 911 & $\lambda=5\cdot 10^{-4}$ \\
&SGP + TV & 47.58 & 841 & $\rho=0.01$ 
& 48.17 & 2526 & $\rho=0.1$
& 56.22 & 3998 & $\rho=0.1$ \\
&SGP & 48.88 & 7000 &  & 37.55 & 546 &  & 44.57 & 140 &  \\
\hline
\multirow{5}{*}{\rotatebox{90}{Rel.\ err.}}&
SGP + Sh + TV  & 0.01 & 1277 & $\lambda=5\cdot 10^{-4}$ 
& 0.10 & 1951 & $\lambda=5\cdot 10^{-4}$
& 0.24 & 432 & $\lambda=5\cdot 10^{-4}$ \\
&& & & $\rho=0.01$ & & & $\rho=1$ & & & $\rho=1$  \\
&SGP + Sh & 0.04 & 426 & $\lambda=5\cdot 10^{-4}$ 
& 0.17 & 745 & $\lambda=5\cdot 10^{-4}$
& 0.91 & 911 & $\lambda=5\cdot 10^{-4}$ \\
&SGP + TV & 0.02 & 841 & $\rho=0.01$ 
& 0.04 & 2526 & $\rho=0.1$
& 0.04 & 3998 & $\rho=0.1$ \\
&SGP & 0.02 & 7000 &  & 0.14 & 546 &  & 0.16 & 140 &  \\
\hline
\end{tabular}}
\caption{Optimal results for the \emph{implicit} formulation of the objective function. The corresponding reconstructed images are reported in Figure 
\ref{fig:recIm}.
}
\label{tab:implicit}
\end{table}

In this section, we present extensive numerical results for the ROI CT 
reconstruction problem in the framework of 2D fan-beam geometry.
All the algorithms were implemented in Matlab 8.1.0 and the experiments 
performed on a dual CPU server, equipped with two 6-cores Intel 
Xeon X5690 at 3.46GHz, 188 GB DDR3 central RAM memory 
and up to 12 TB of disk storage. 
The object to be imaged is the modified Shepp-Logan phantom sized 
$N \times N$ pixels with $N=128$. We assumed it to be placed in the first 
quadrant of the Cartesian coordinate system with the image lower left corner 
at the origin. 
Truncated projection data were obtained by discarding the samples 
outside the ROI projection $\mcP(S)$. We recall that this corresponds to a 
ROI disk in the image domain (see Section \ref{sec:contsetup}). In particular, 
we considered concentric ROI disks with decreasing \textit{radius} 
$\gamma$, placed off-center with respect to the field of view. 
The results in this paper covers ROIs that are fully inside the 
object being imaged (in the hypothesis of the interior tomography problem). 
However, a larger \textit{radius} that exceeds the object along one 
coordinate axis has been considered for comparison.  

\begin{figure}
\vspace*{-7em}
\centering
\makebox[0pt][c]{%
\begin{tabular}{@{}ccc@{}}
$\gamma = 0.5 N$ & $\gamma = 0.3 N$ & $\gamma = 0.15 N$  \\
\includegraphics[scale=0.33]{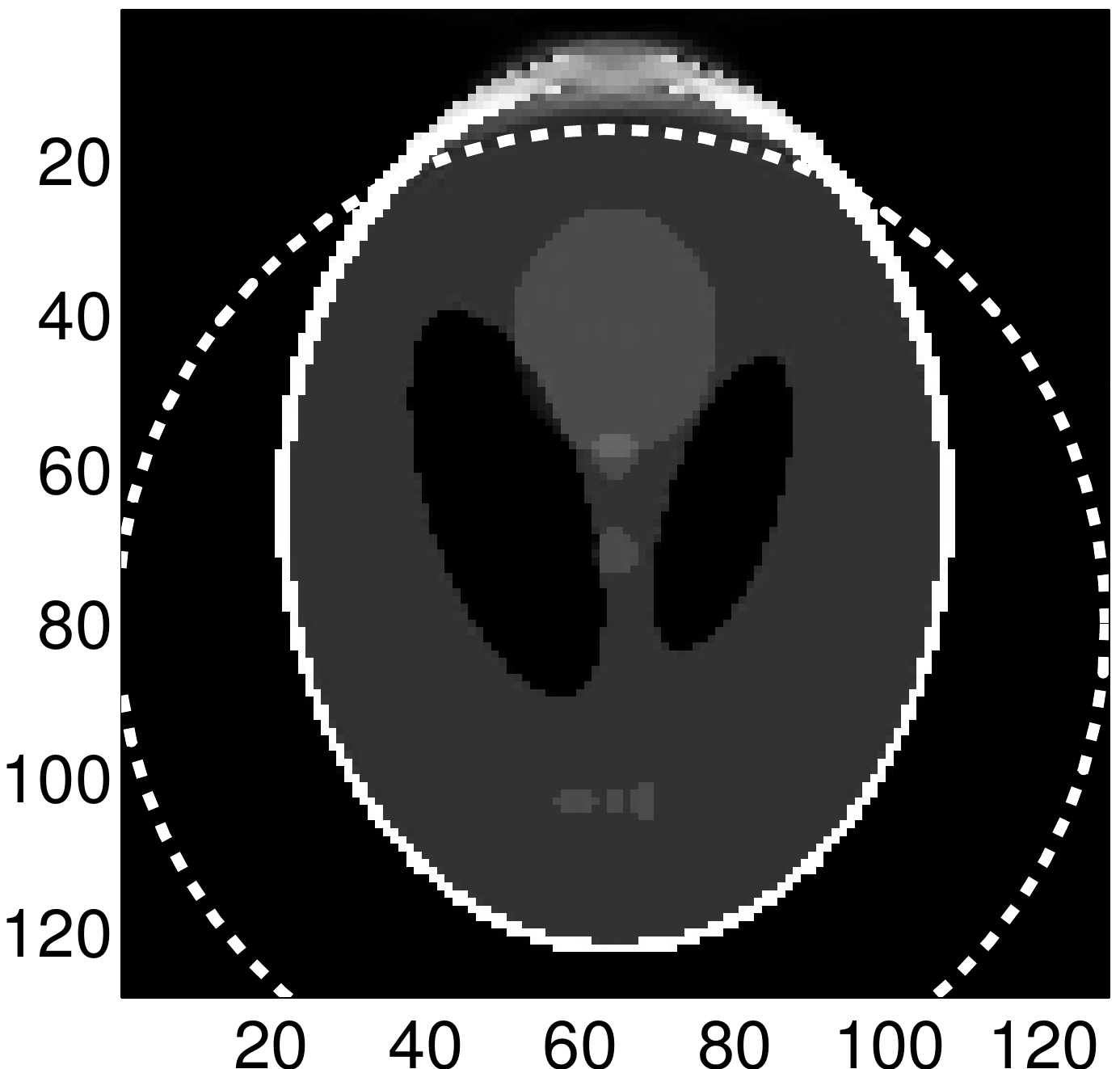}
& \includegraphics[scale=0.33]{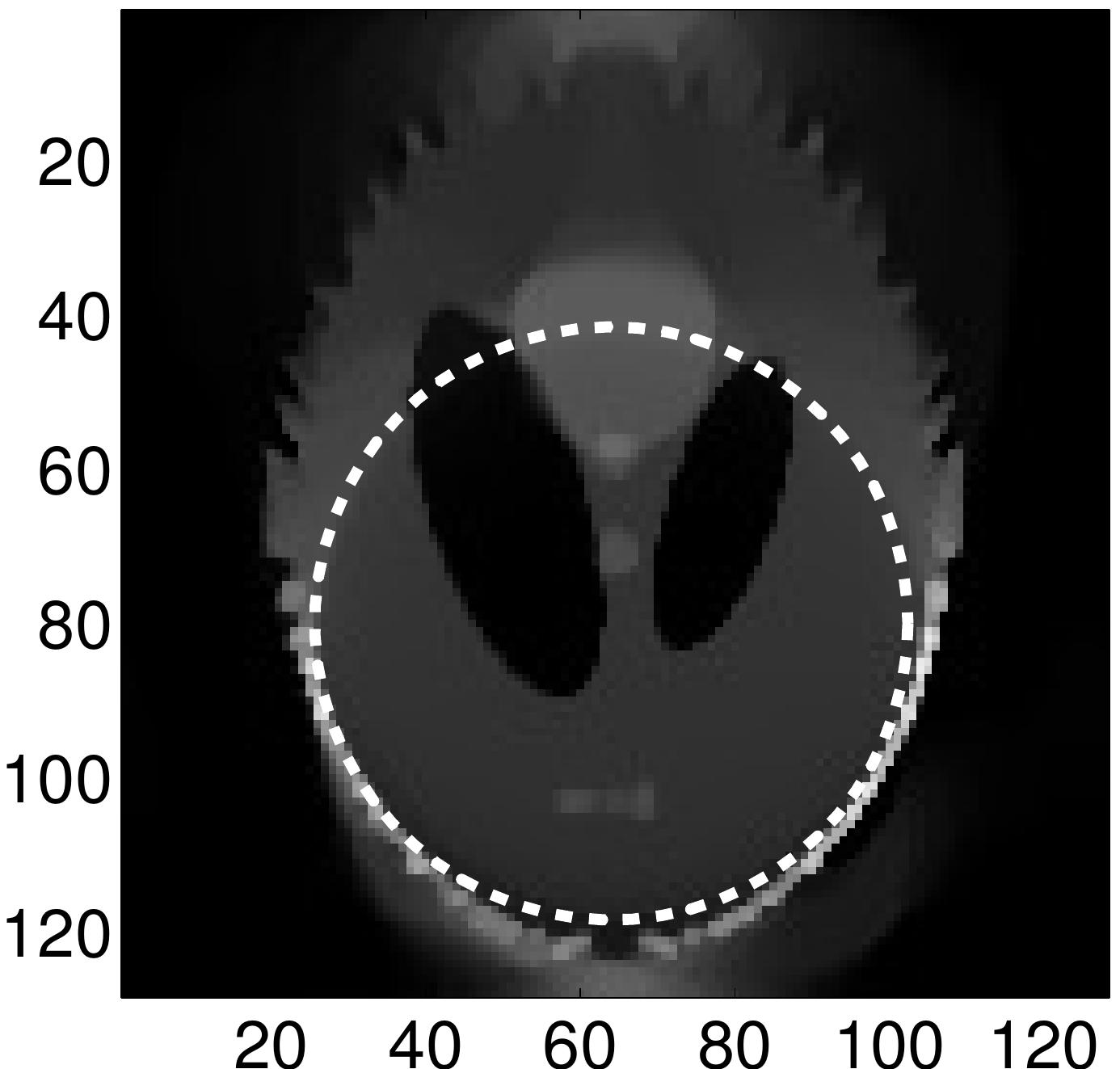}
& \includegraphics[scale=0.33]{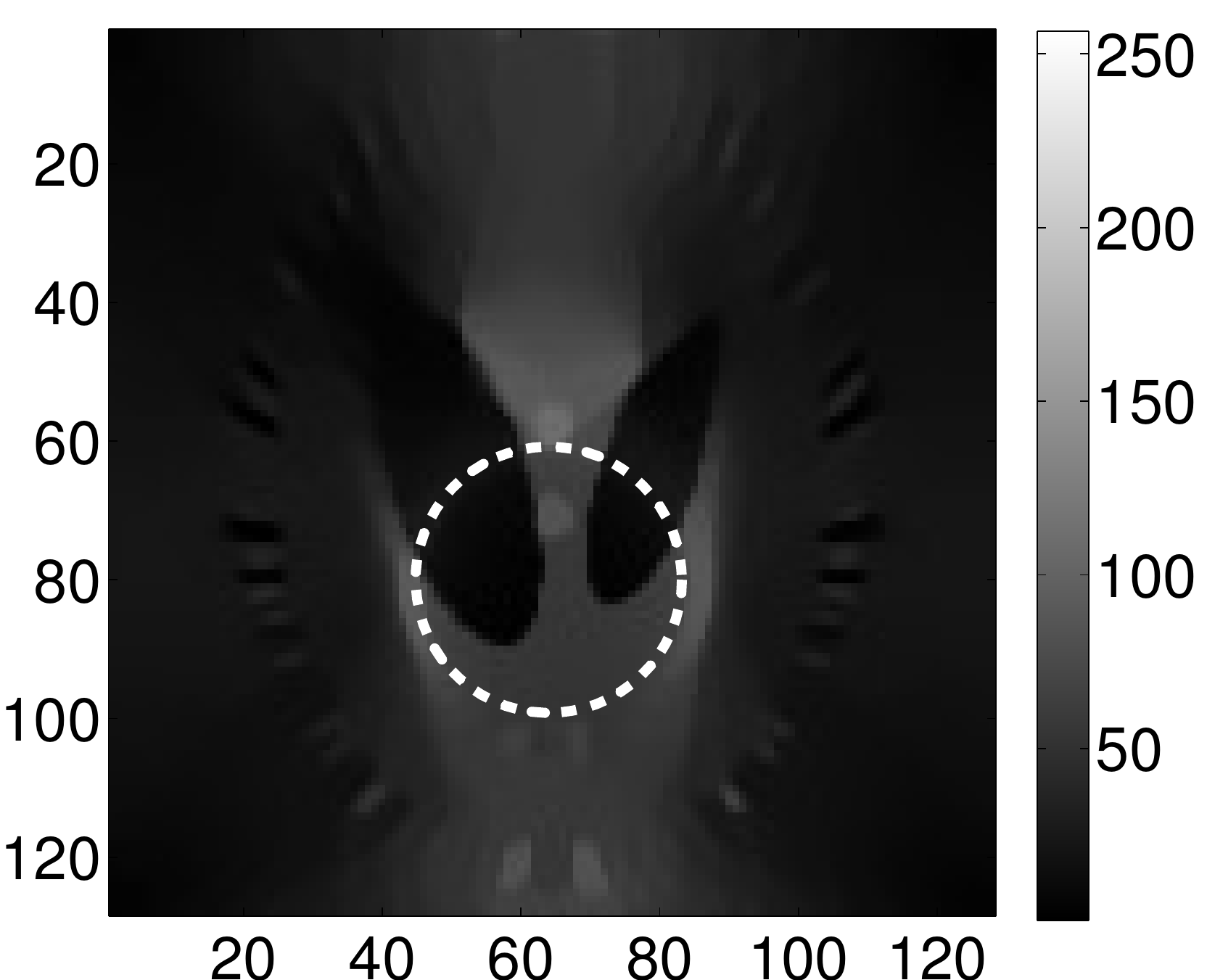} \\
(a) & (b) & (c) \\
\includegraphics[scale=0.33]{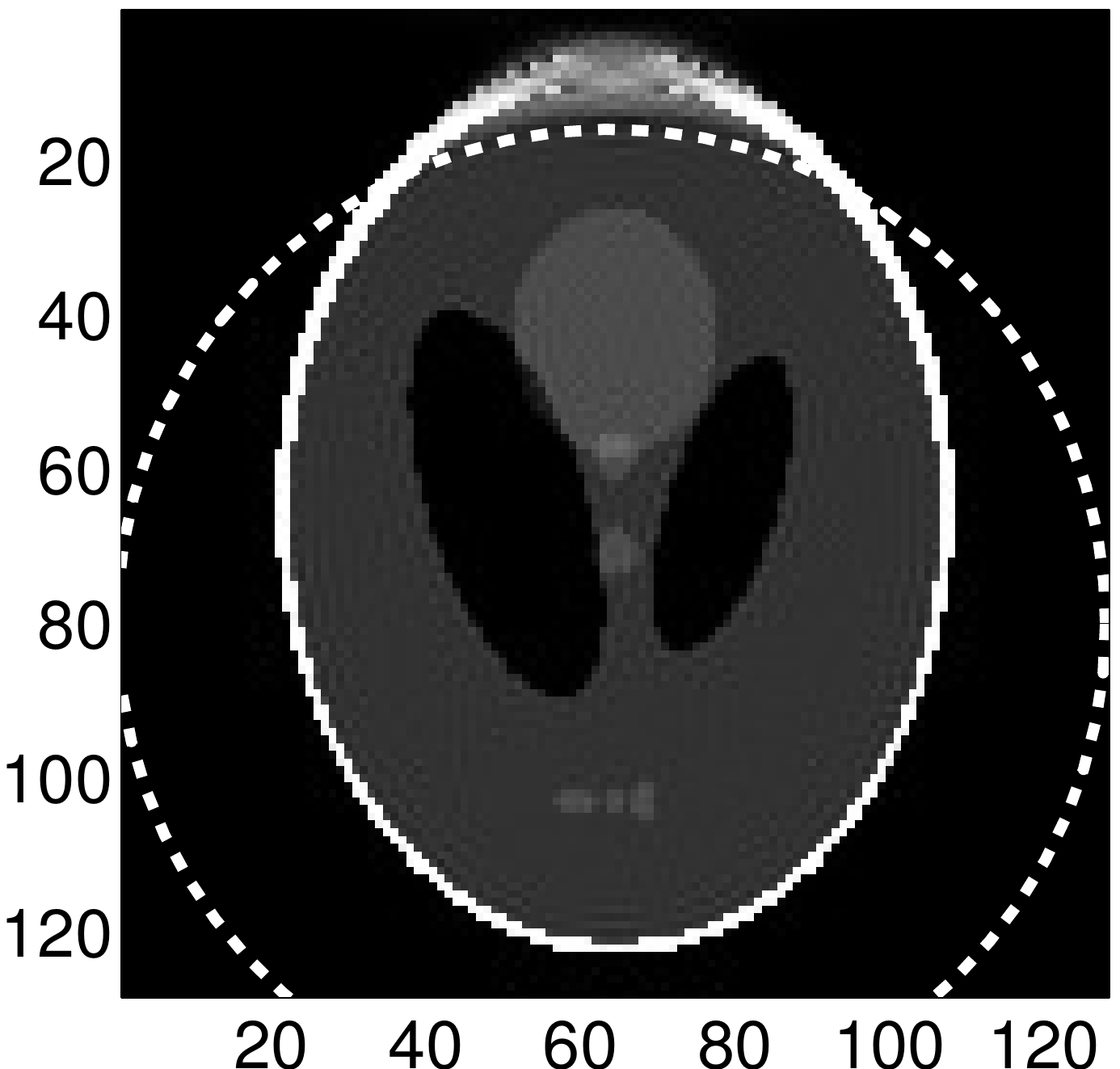}
& \includegraphics[scale=0.33]{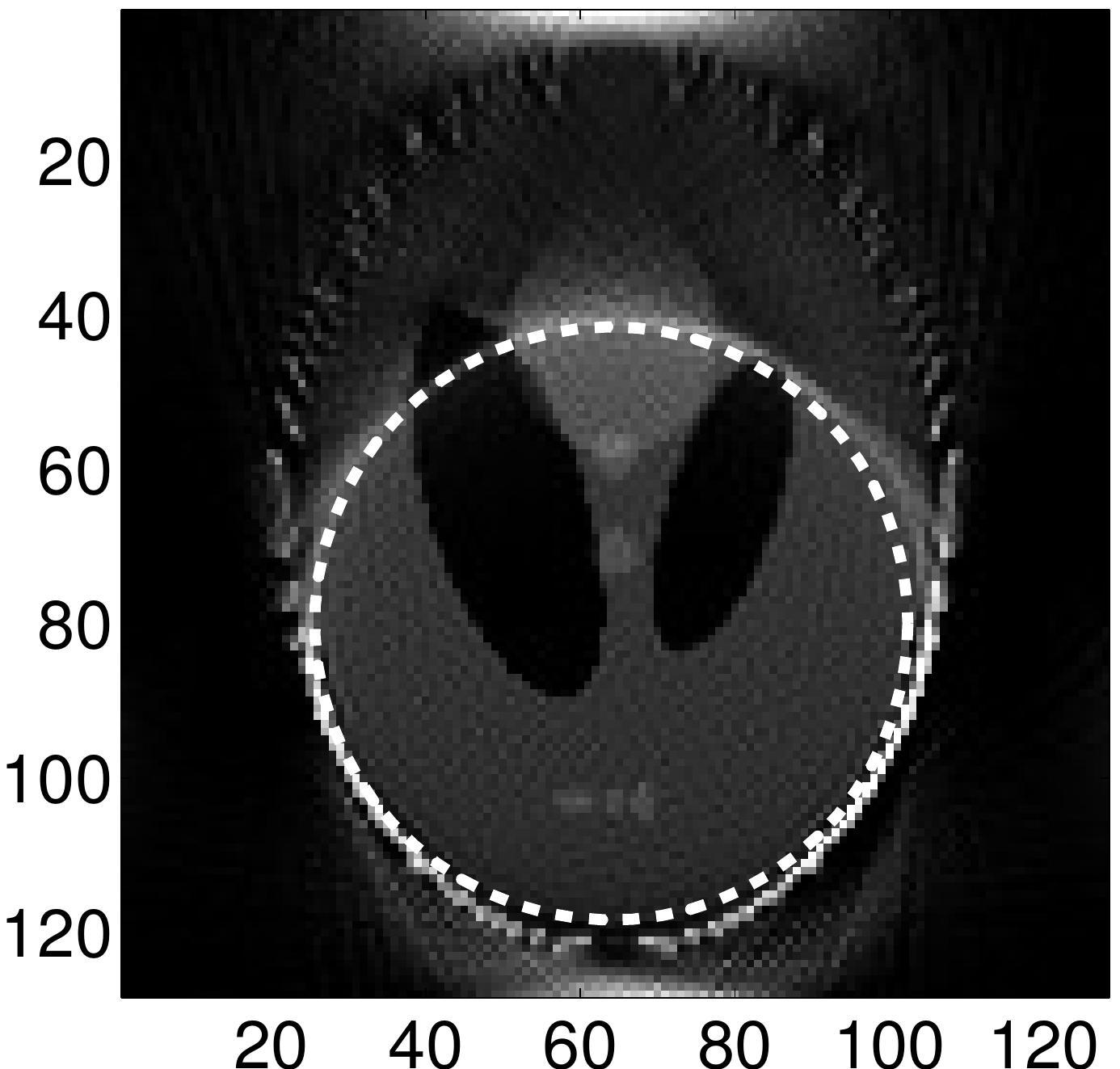}
& \includegraphics[scale=0.33]{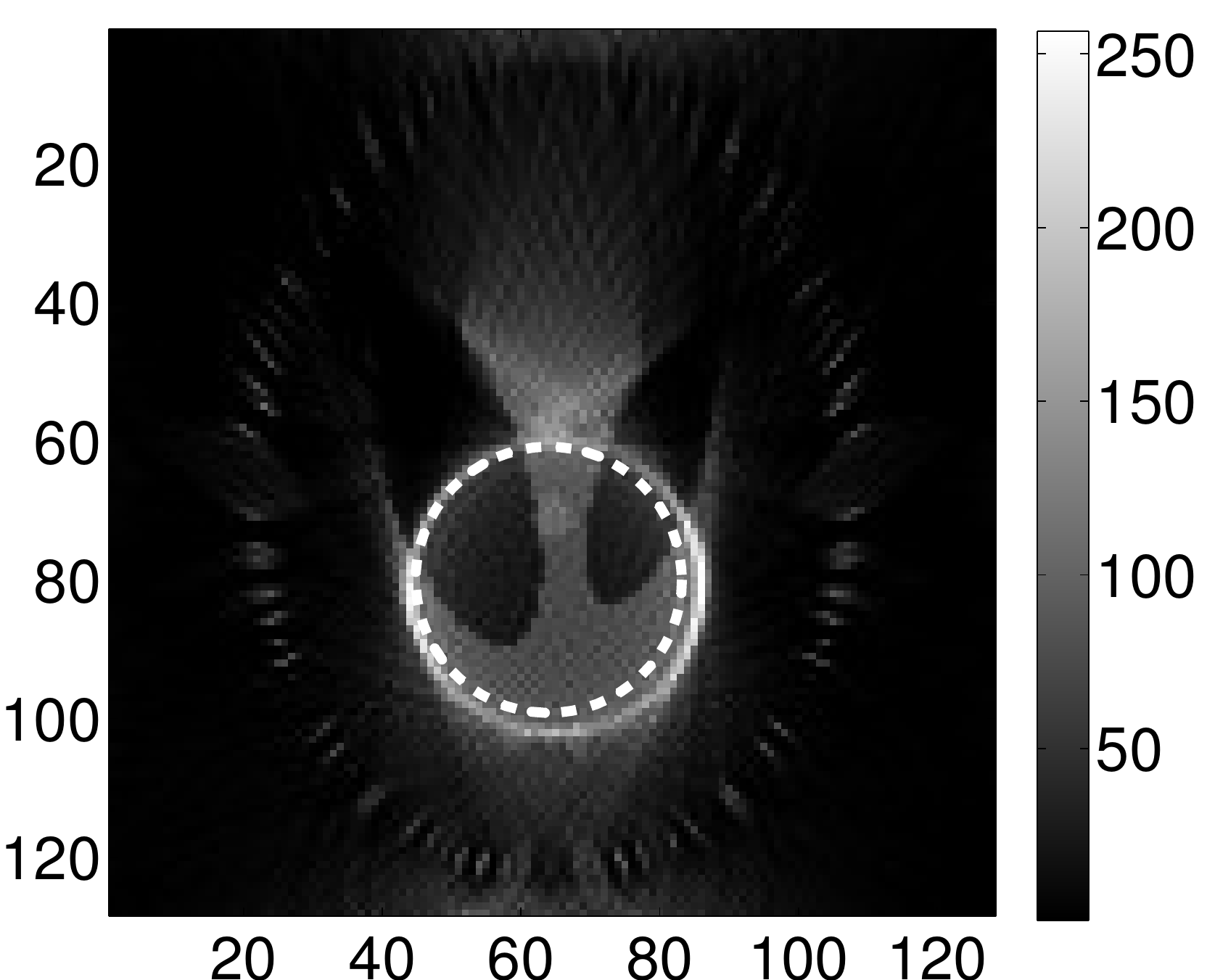} \\
(d) & (e) & (f)  \\
\includegraphics[scale=0.33]{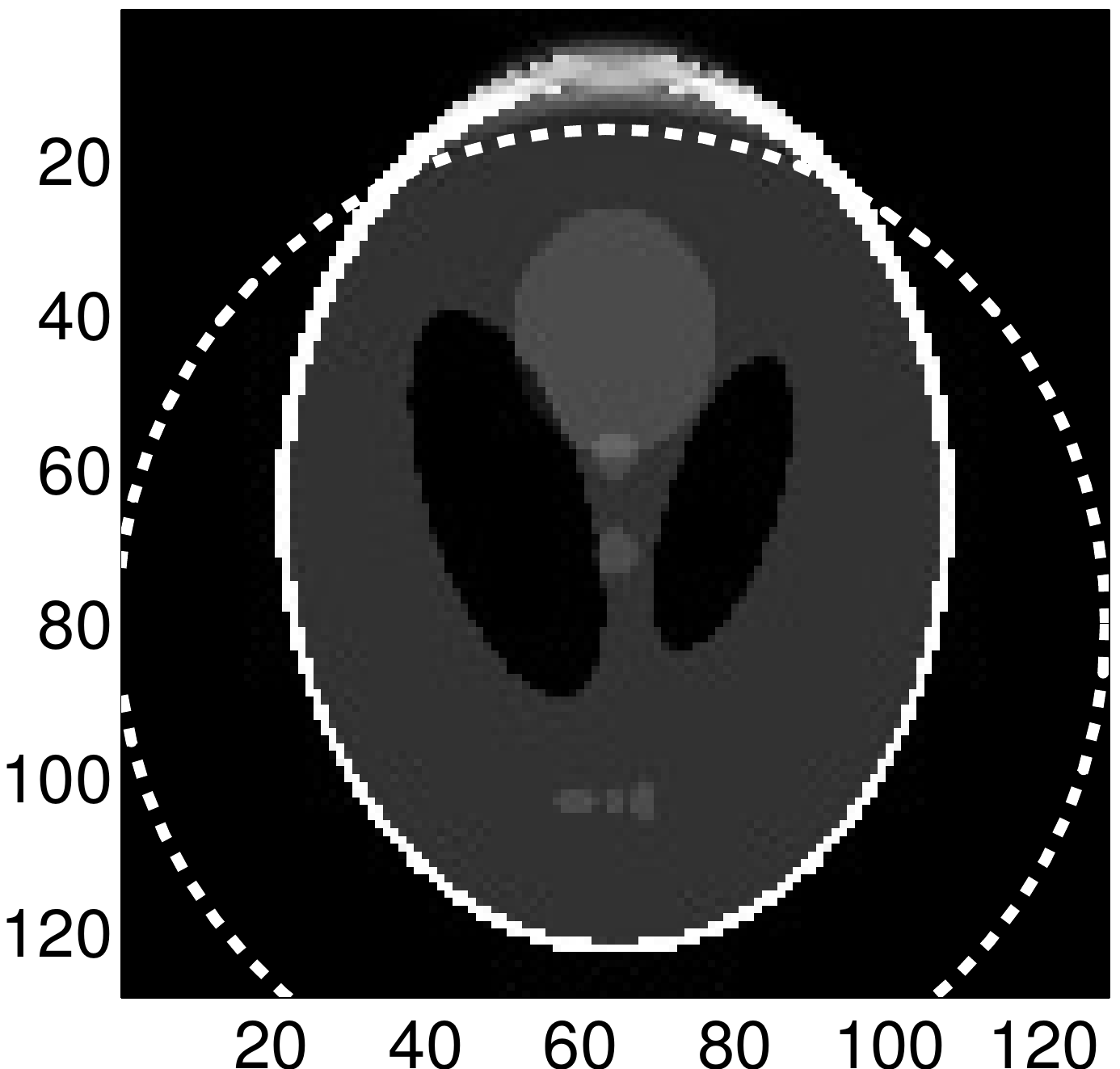}
& \includegraphics[scale=0.33]{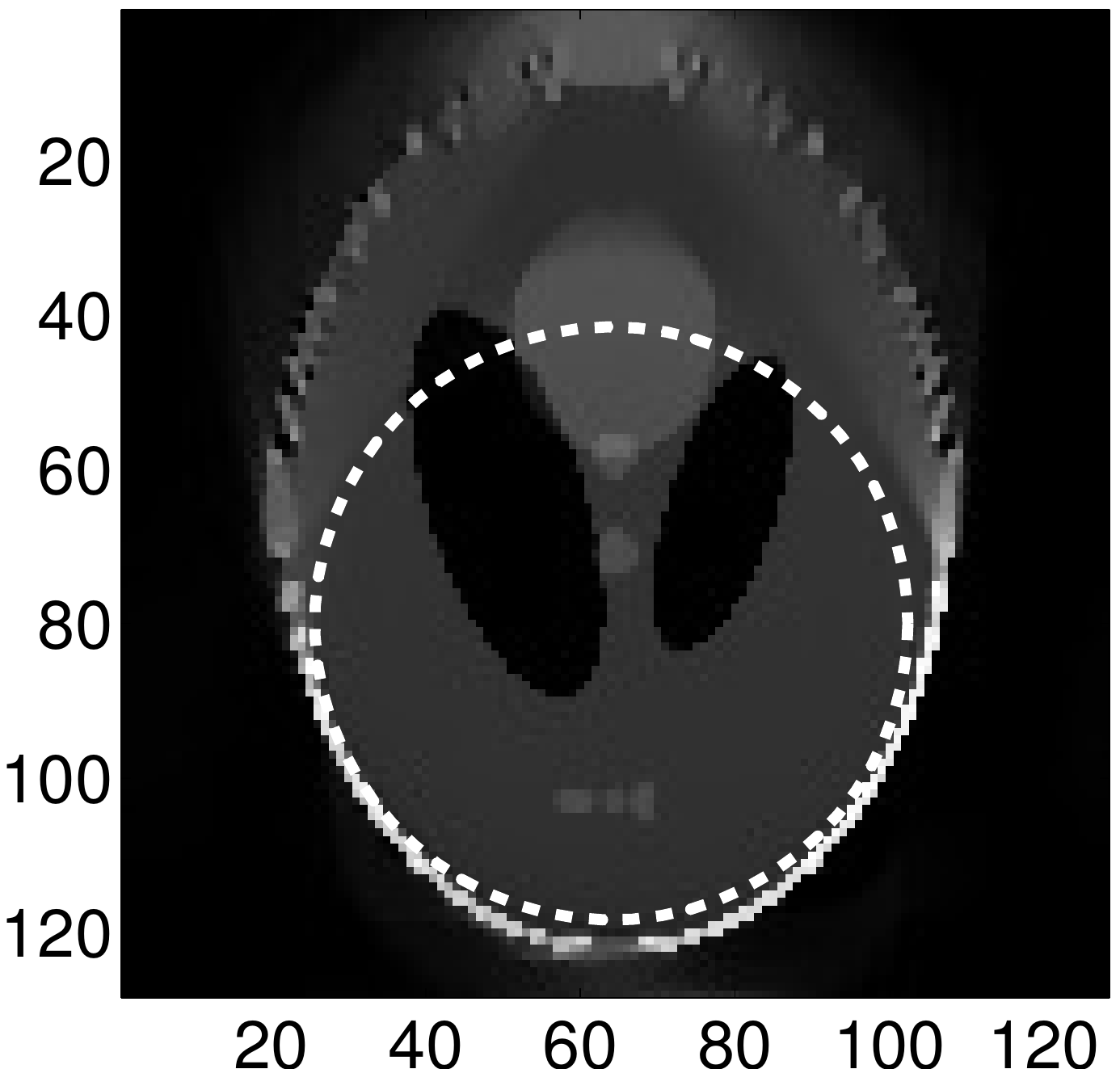}
& \includegraphics[scale=0.33]{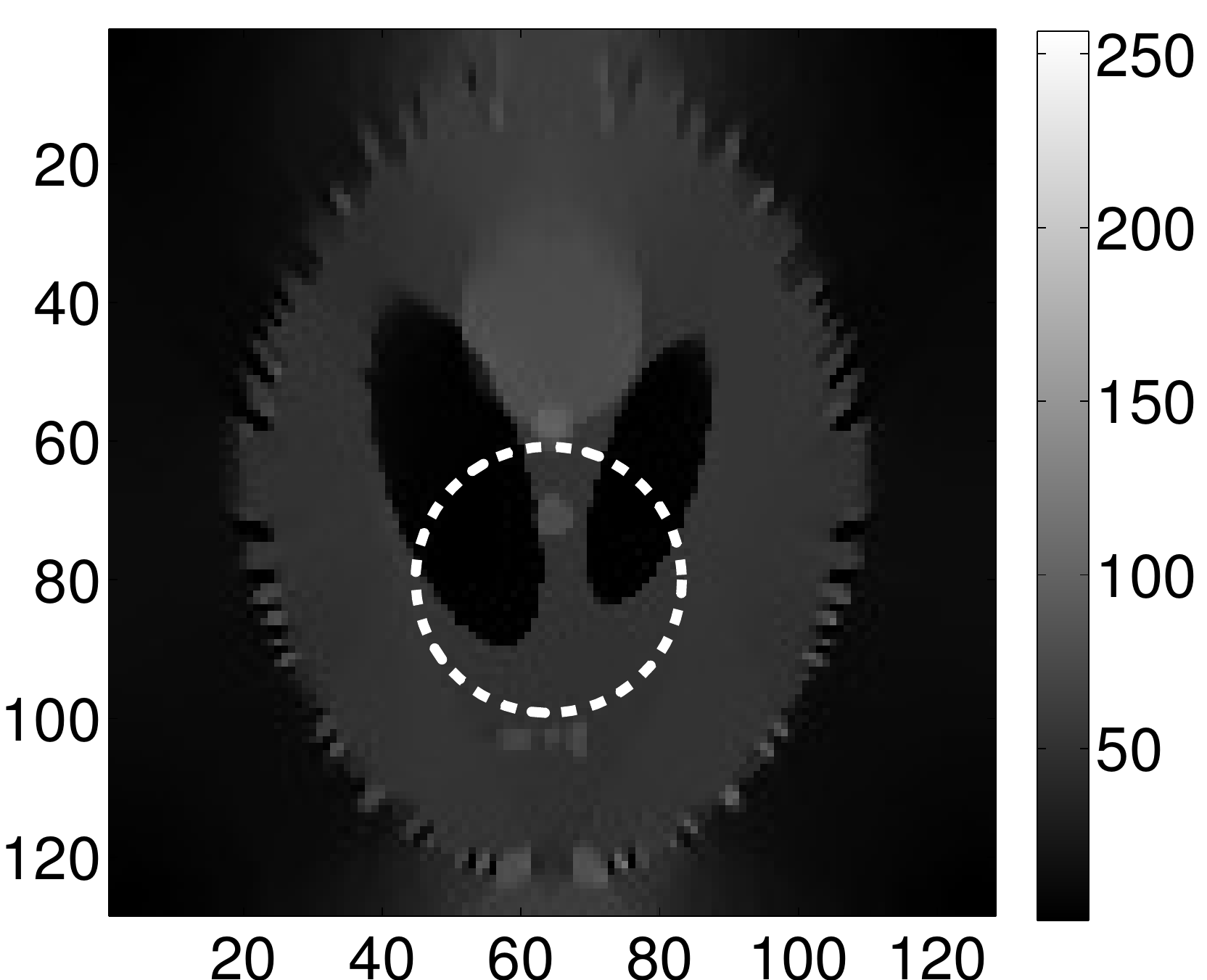} \\
(g) & (h) & (i)   \\
\includegraphics[scale=0.33]{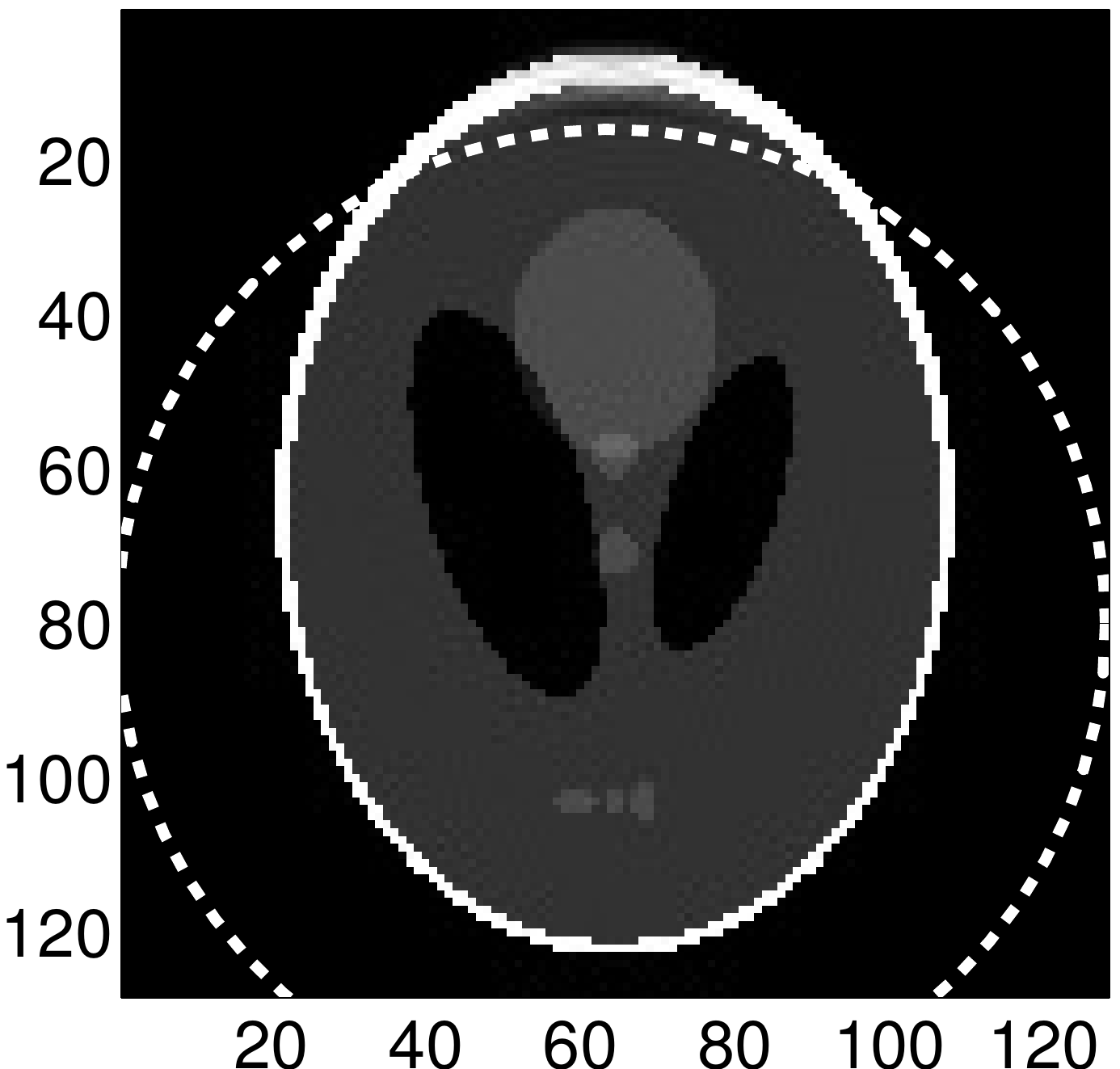}
& \includegraphics[scale=0.33]{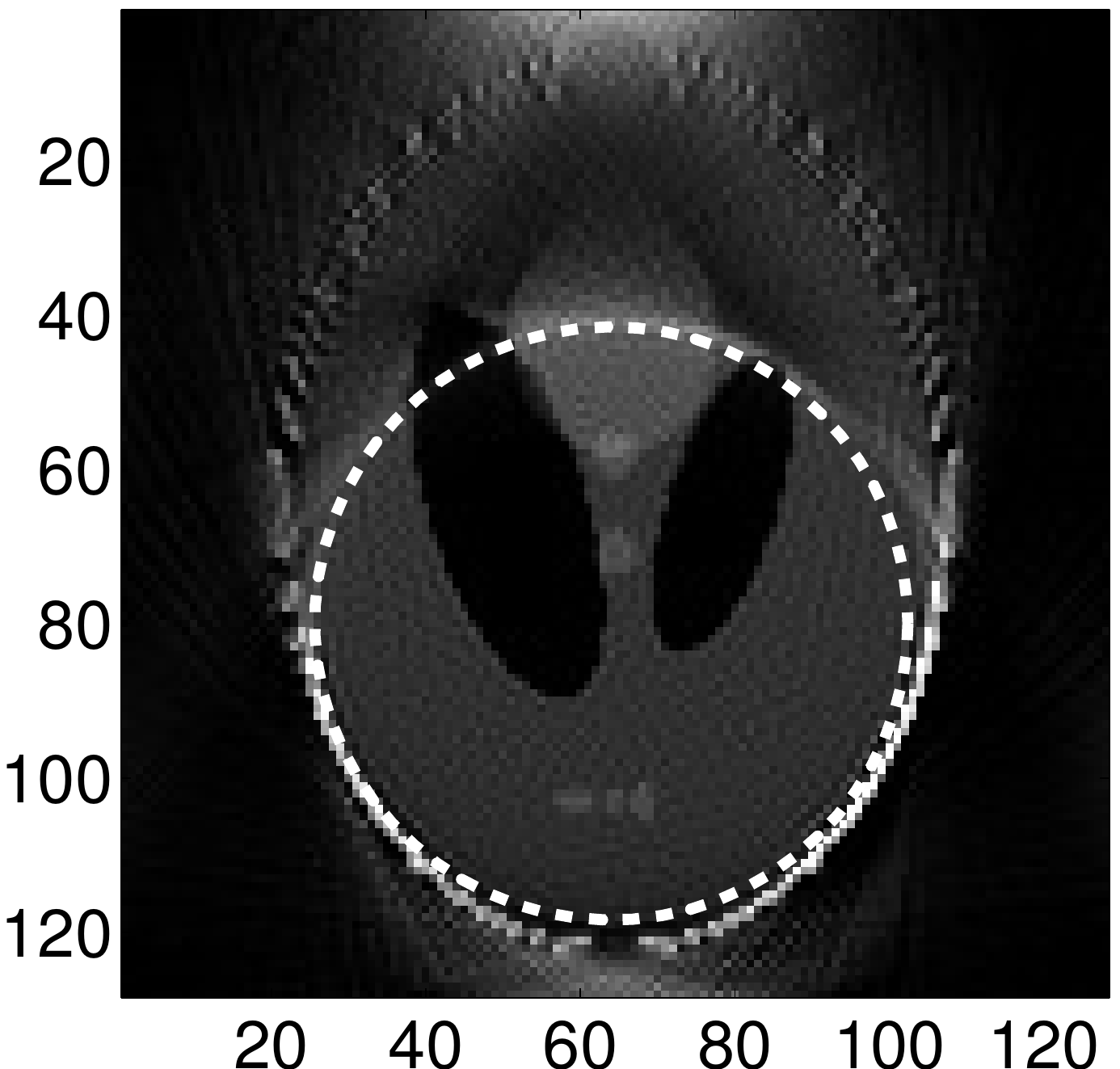}
& \includegraphics[scale=0.33]{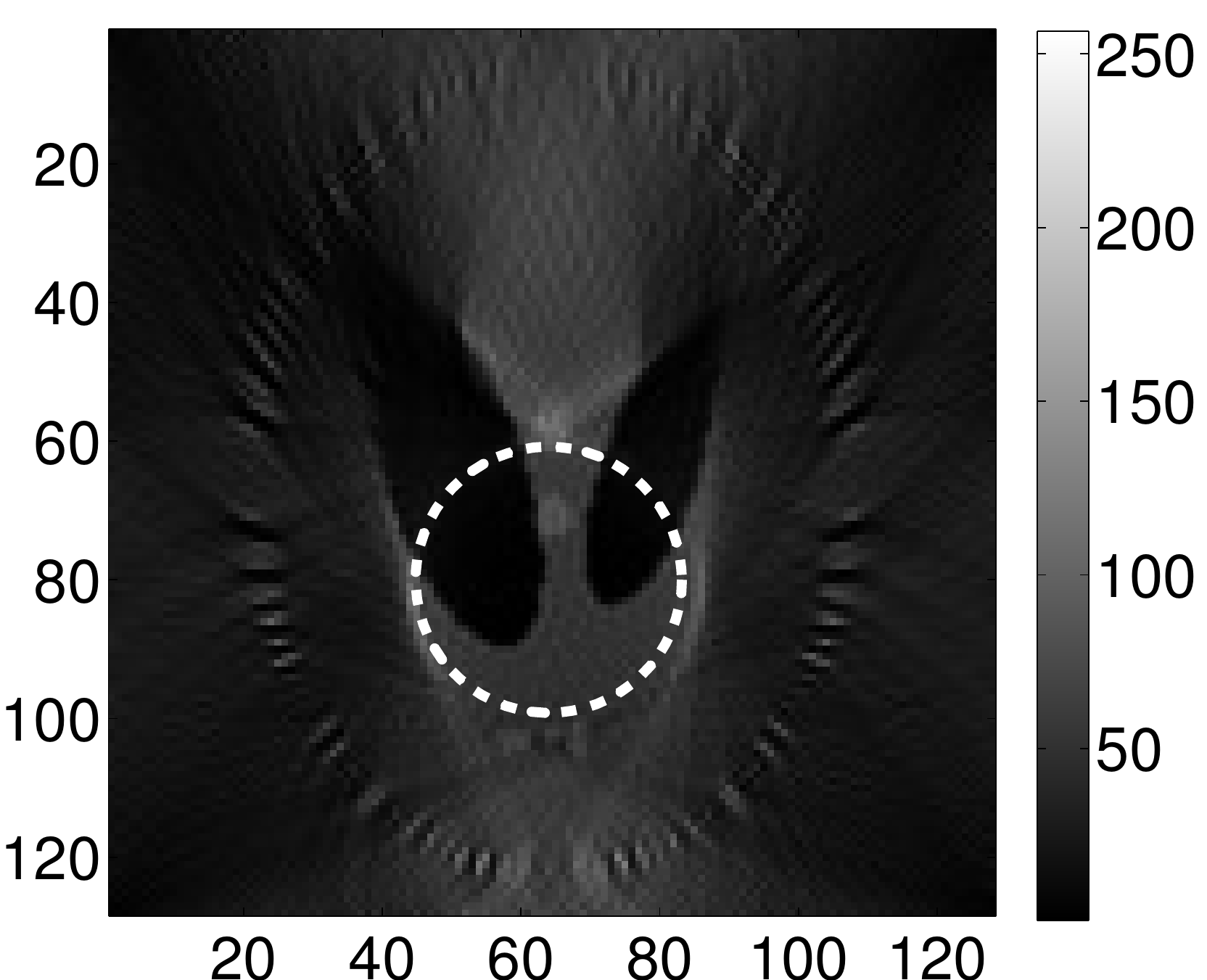} \\
(j) & (k) & (l)
\end{tabular}}
\caption{Optimal reconstructions of the Shepp-Logan phantom for implicit formulation with decreasing 
\textit{radii}: $\gamma=0.5 N$ for (a), (d), (g) and (j), $\gamma=0.3 N$ for (b), (e), (h) and (k), $\gamma=0.15 N$ for (c), (f), (i) and (l). 
First row: shealerts and TV. Second row: just shearlets. 
Third row: just TV. Fourth row: early stopping only. 
}
\label{fig:recIm} 
\enlargethispage{6\baselineskip}
\end{figure}

\begin{figure}
\vspace*{-7em}
\centering
\makebox[0pt][c]{%
\begin{tabular}{@{}cccc@{}}
$\gamma = 0.5 N$ & $\gamma = 0.3 N$ & $\gamma = 0.15 N$  \\
\includegraphics[scale=0.33]{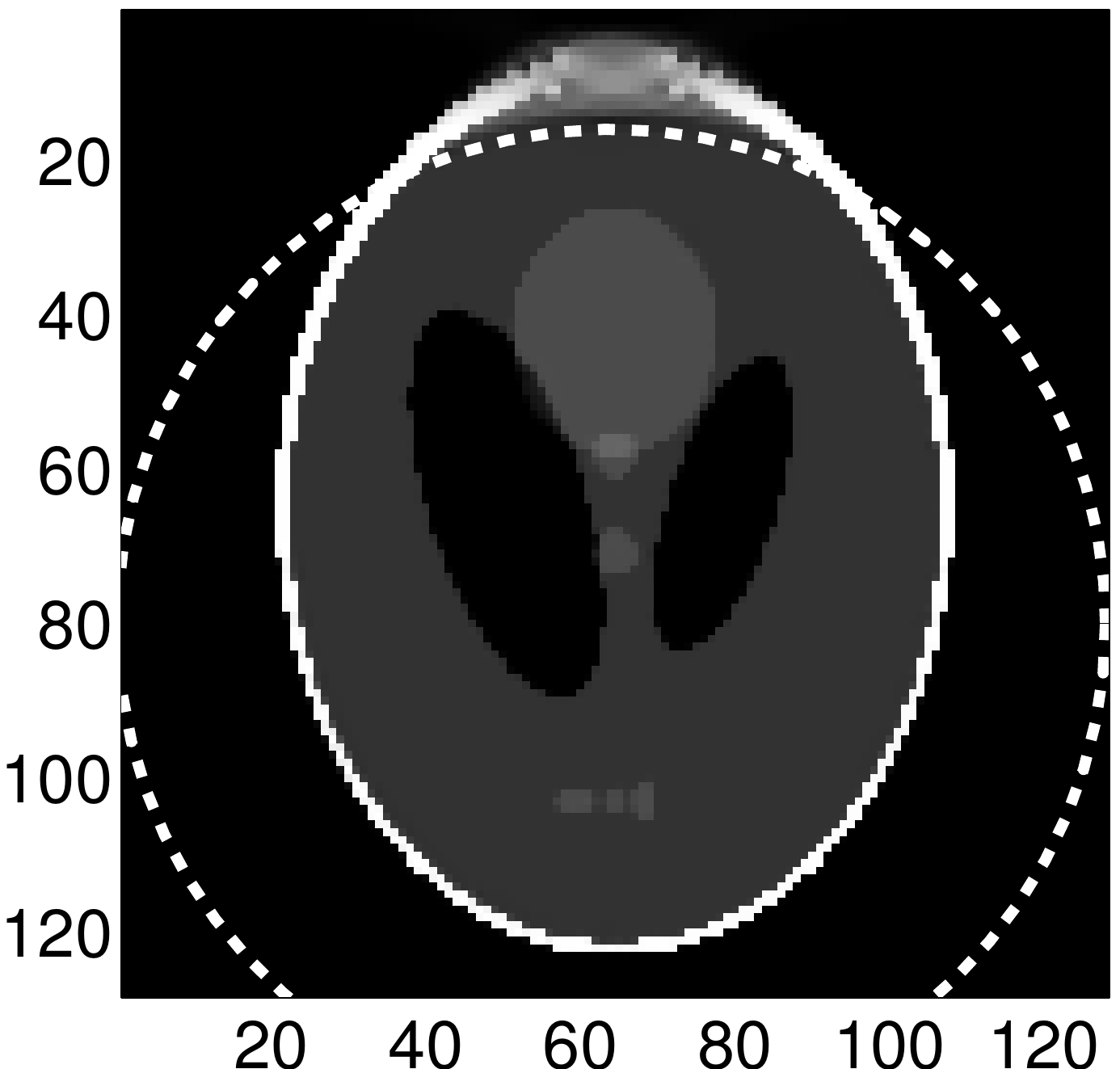}
& \includegraphics[scale=0.33]{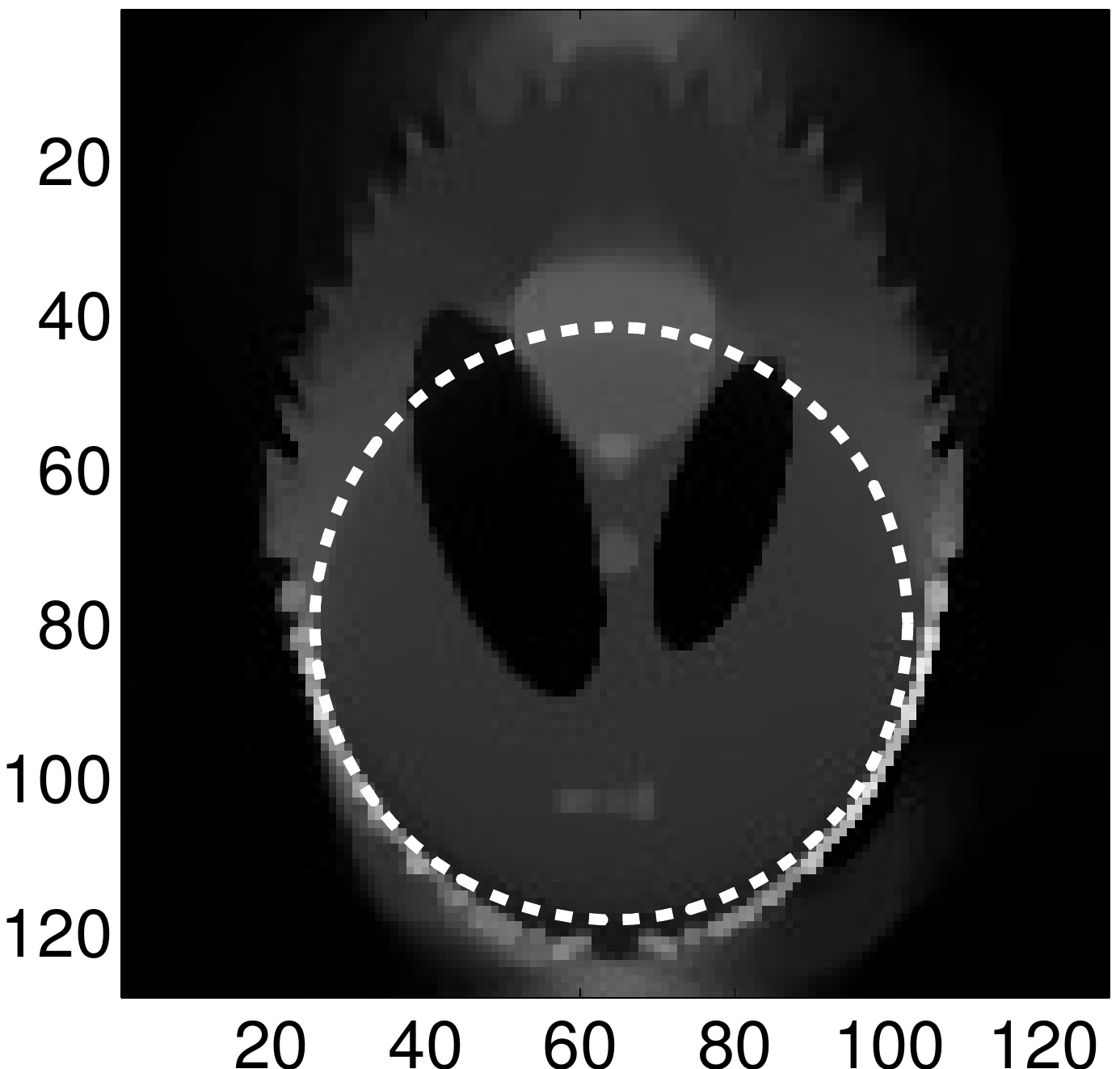}
& \includegraphics[scale=0.33]{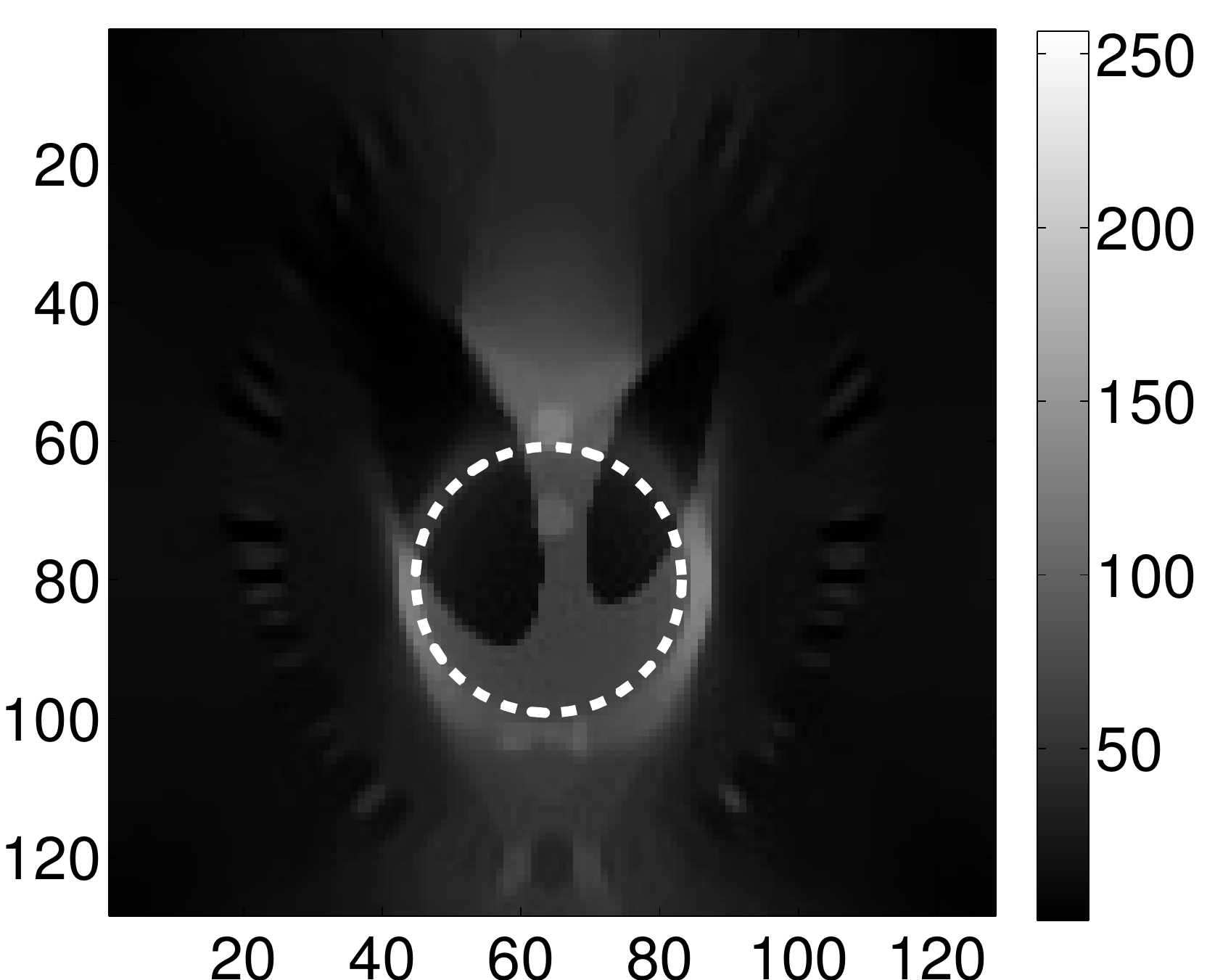} \\
(a) & (b) & (c) \\
\includegraphics[scale=0.33]{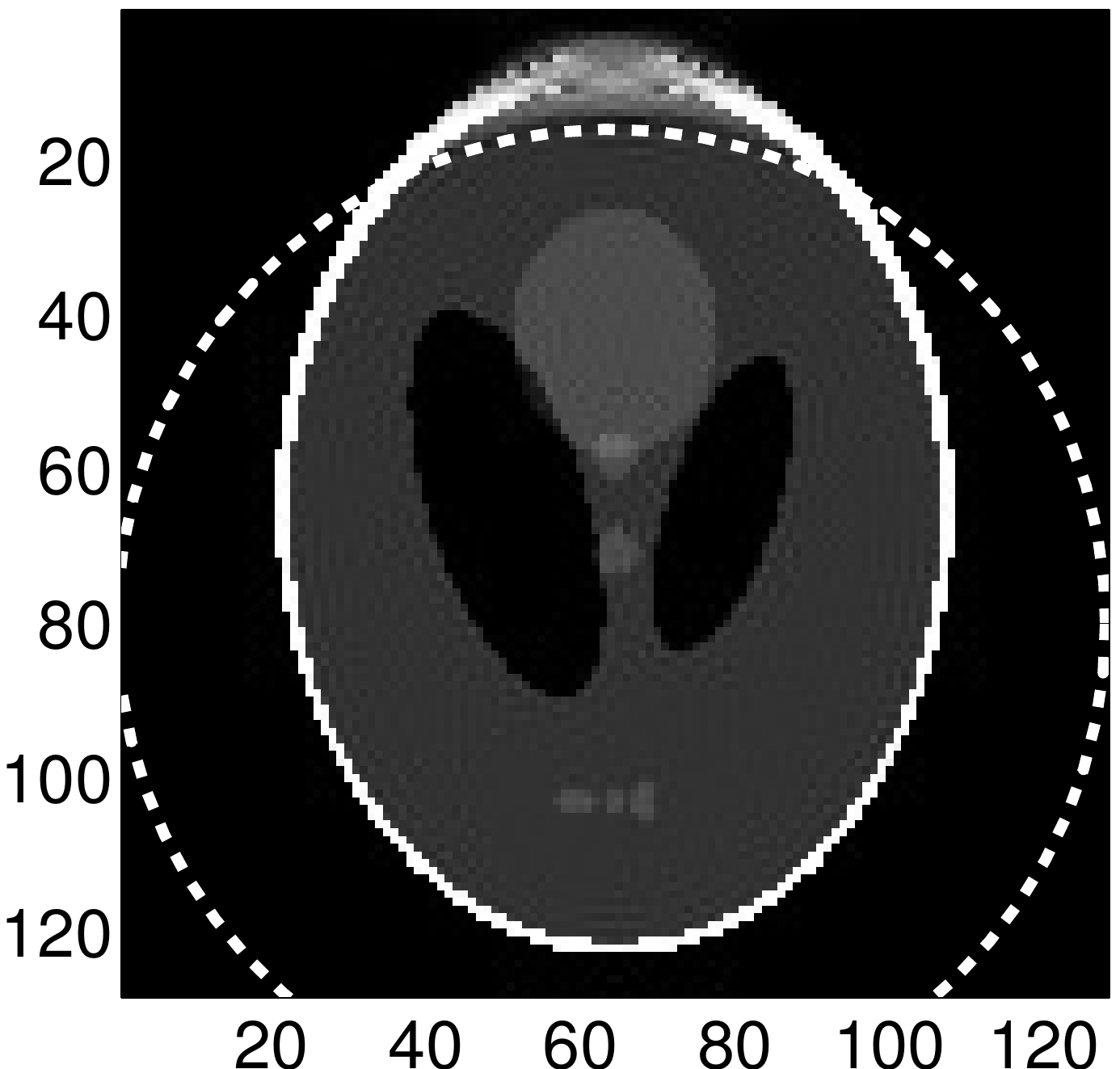}
& \includegraphics[scale=0.33]{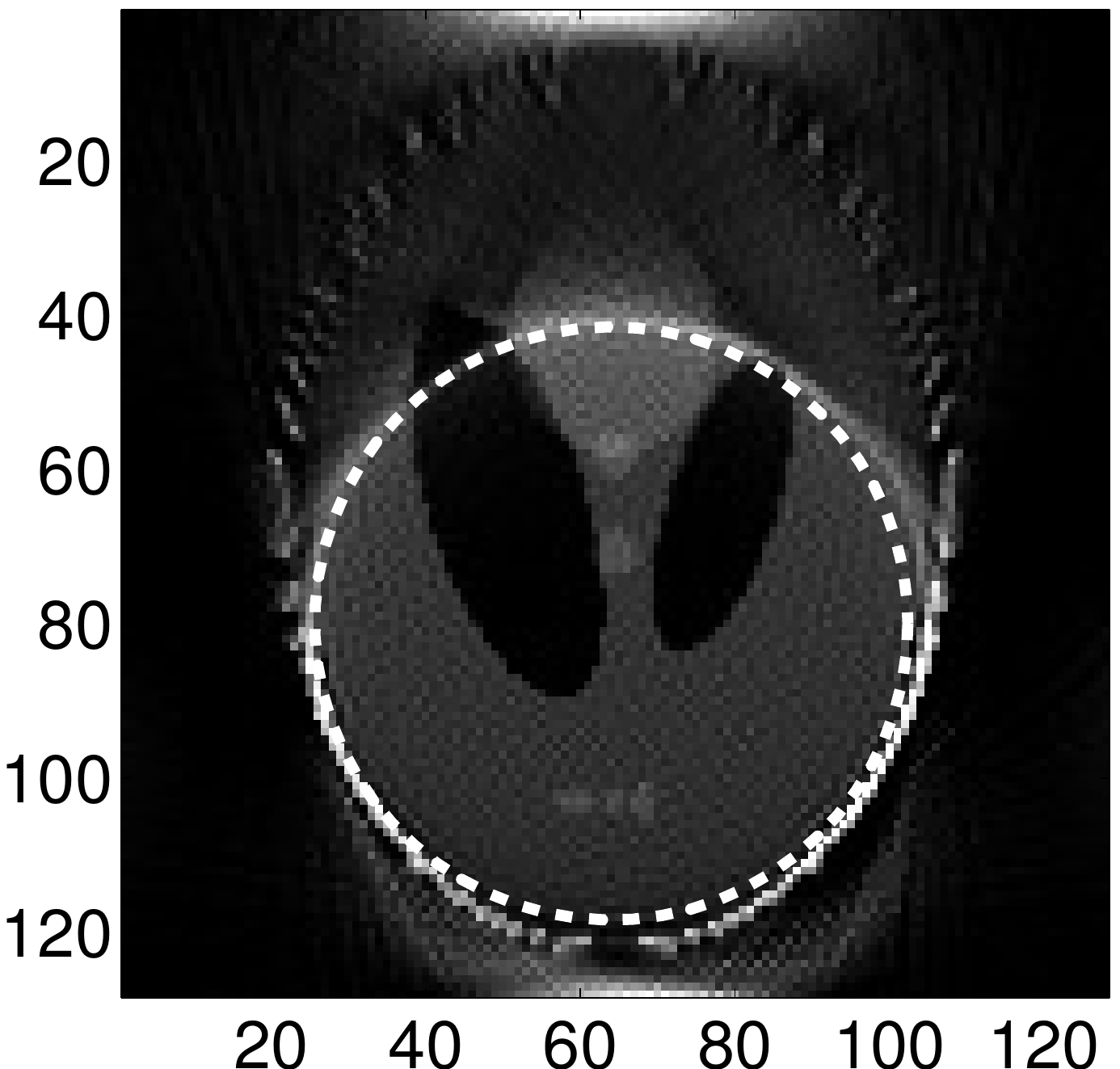}
& \includegraphics[scale=0.33]{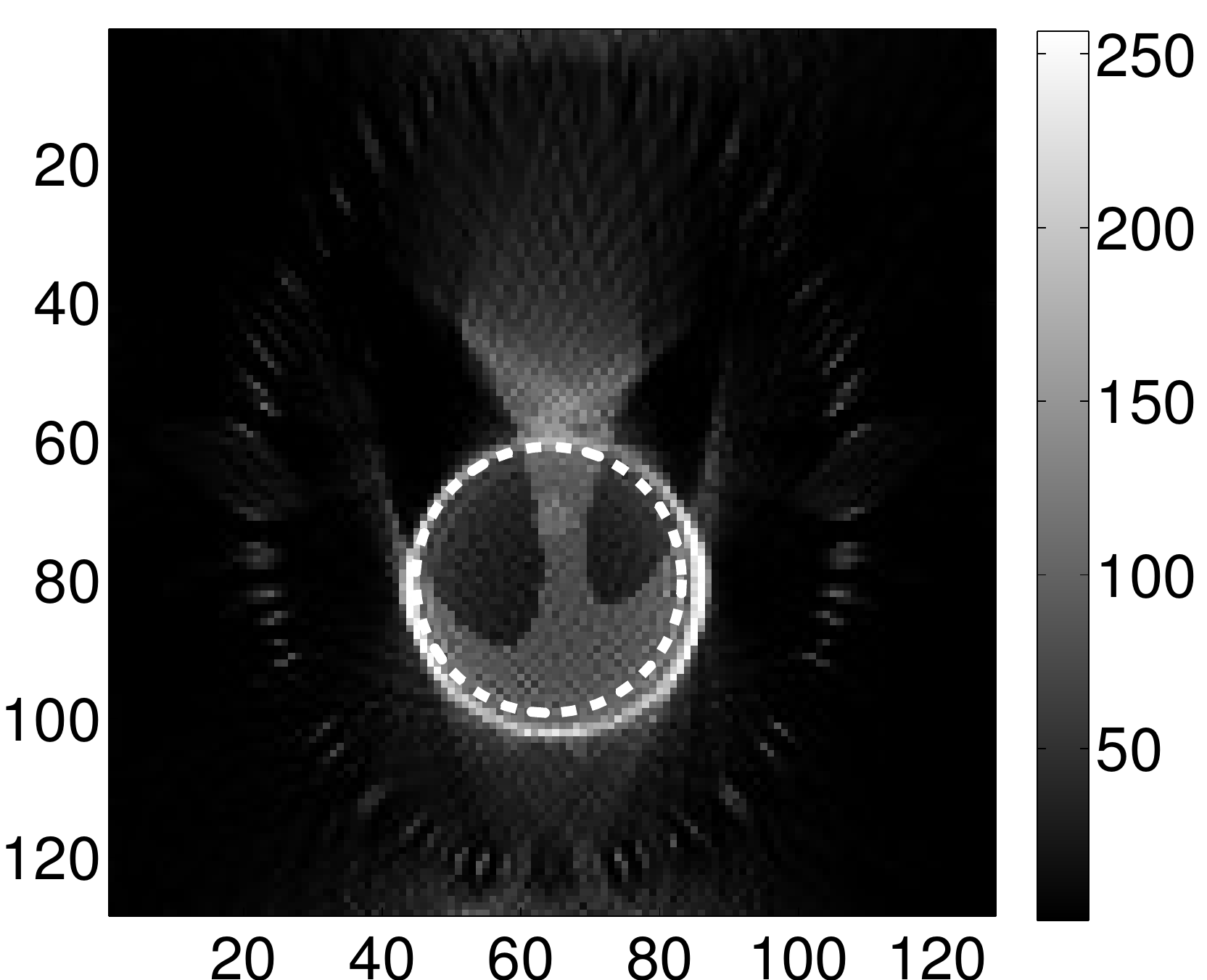} \\
(d) & (e) & (f)  \\
\includegraphics[scale=0.33]{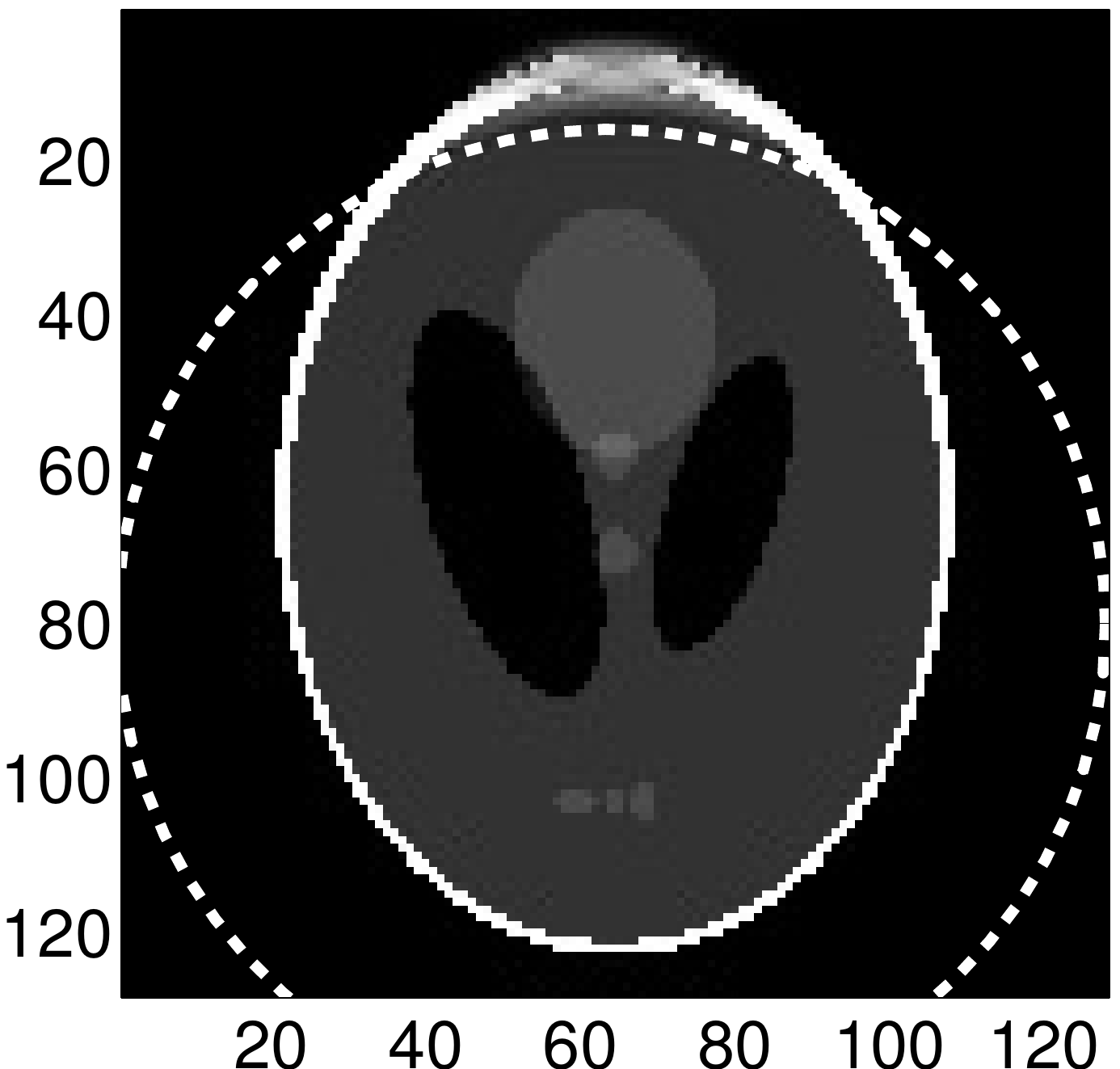}
& \includegraphics[scale=0.33]{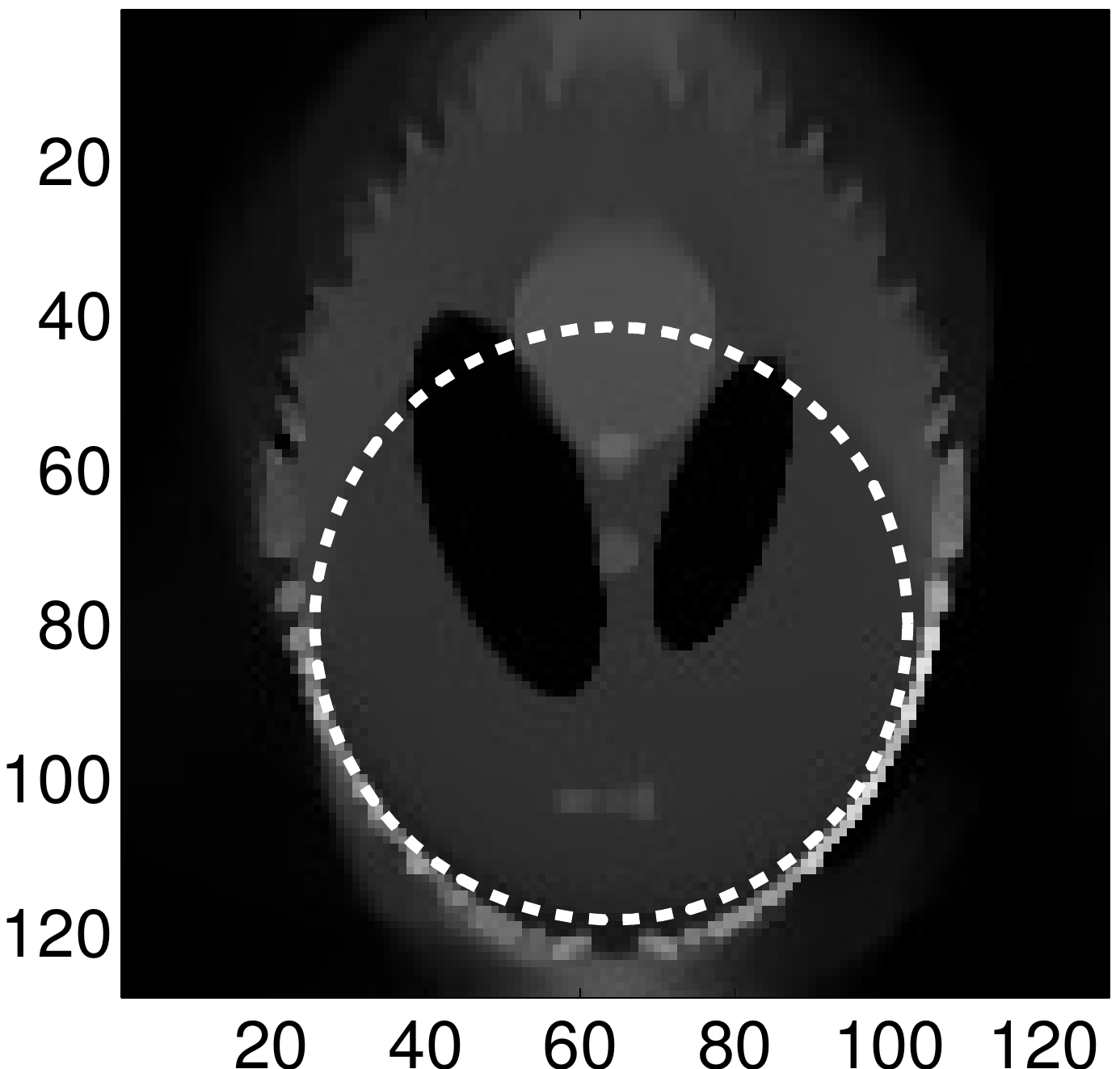}
& \includegraphics[scale=0.33]{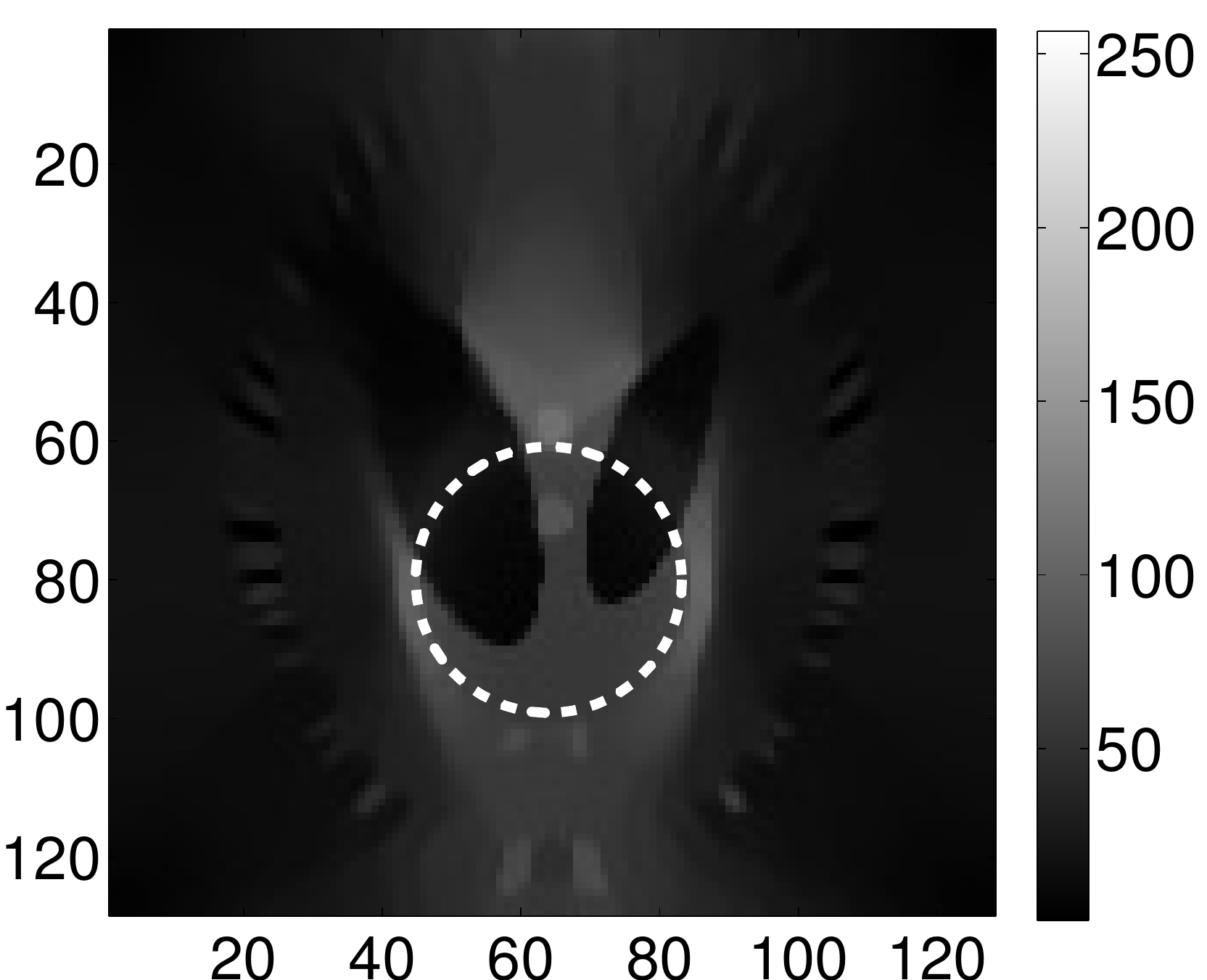} \\
(g) & (h) & (i)   \\
\includegraphics[scale=0.33]{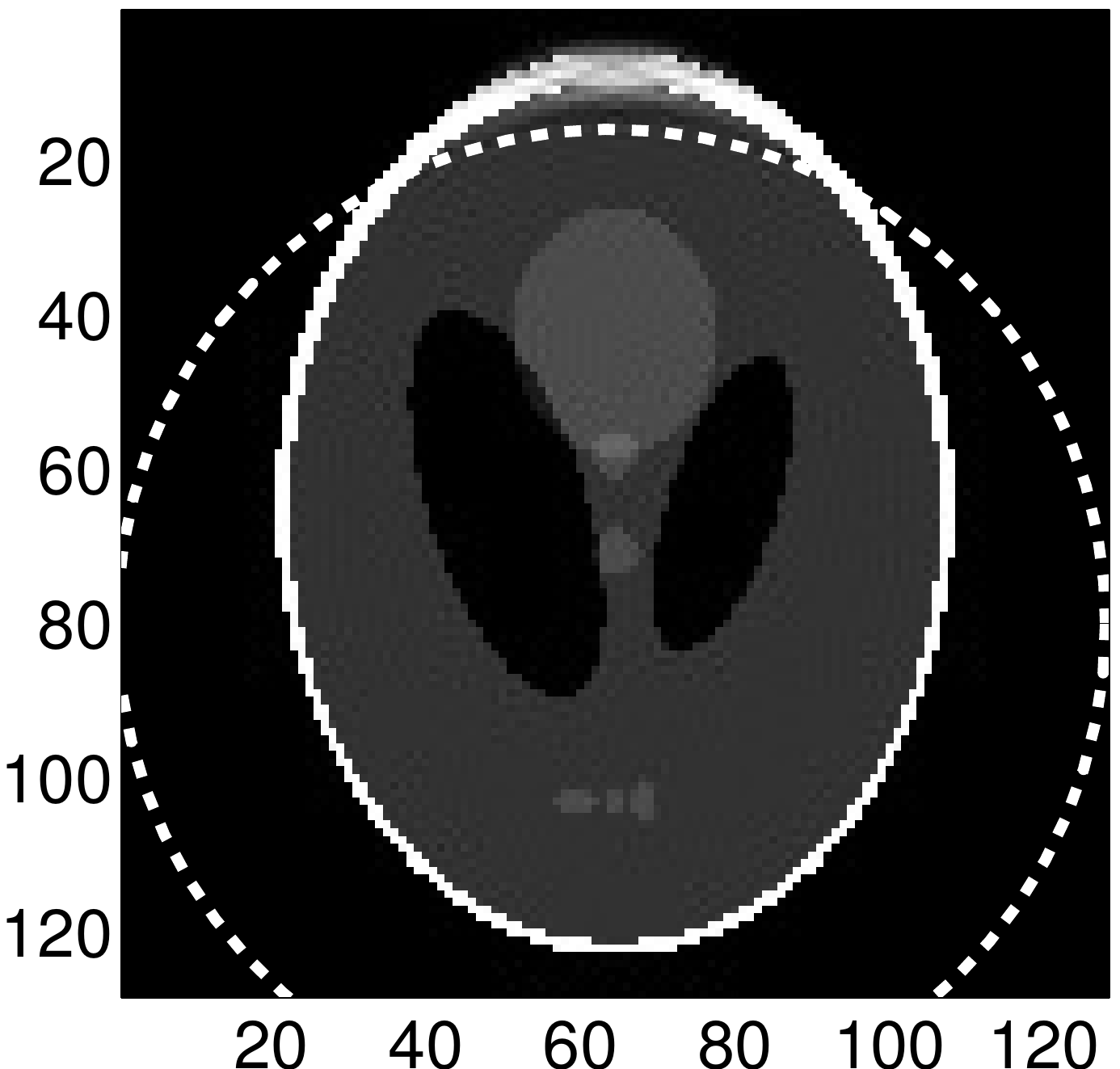}
& \includegraphics[scale=0.33]{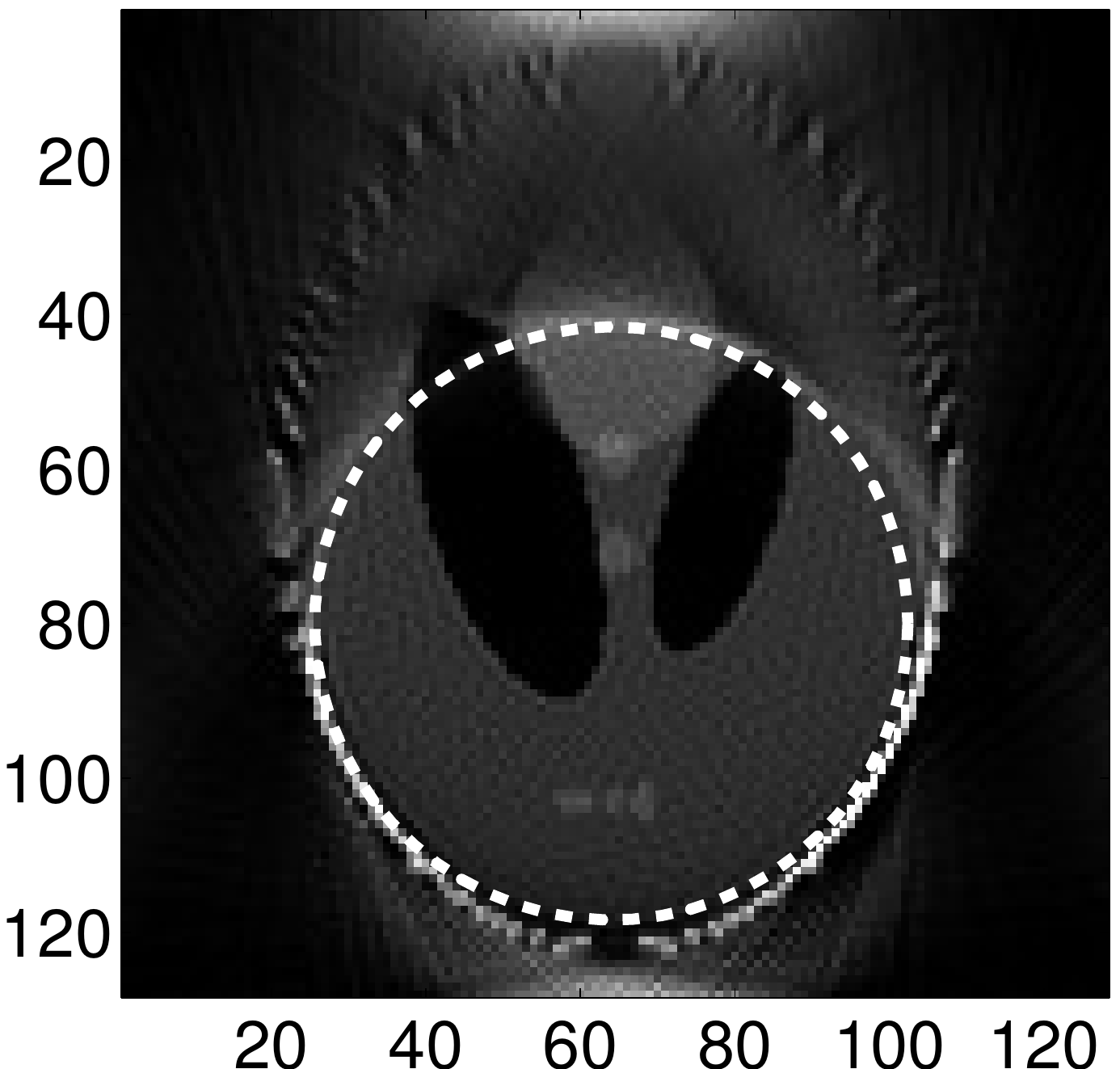}
& \includegraphics[scale=0.33]{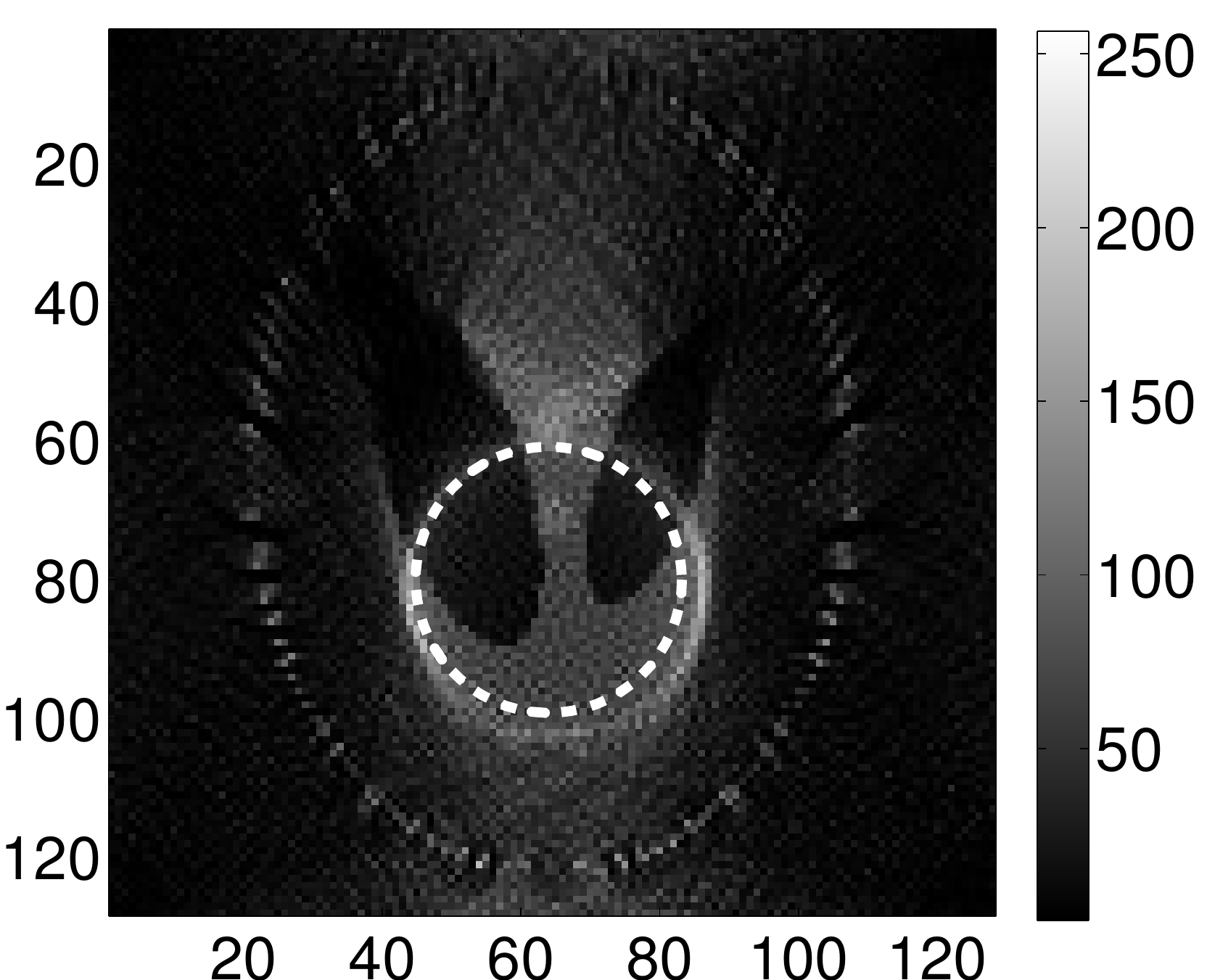} \\
(j) & (k) & (l)
\end{tabular}}
\caption{Optimal reconstructions of the Shepp-Logan phantom for explicit formulation with decreasing 
\textit{radii}: $\gamma=0.5 N$ for (a), (d), (g) and (j), $\gamma=0.3 N$ for (b), (e), (h) and (k), $\gamma=0.15 N$ for (c), (f), (i) and (l). 
First row: shearlets and TV. Second row: just shearlets. 
Third row: just TV. Fourth row: early stopping only. 
}
\label{fig:recEx} 
\end{figure}

\nopagebreak[3]
Numerical tests were performed to assess the best approach
with respect to both the objective function and the regularization parameters.
Tables \ref{tab:implicit} and \ref{tab:explicit} show some of the results for 
the implicit and the explicit formulations, respectively, compared against two 
state-of-the-art figures of merit, namely the peak-signal-to-noise ratio 
(PSNR) and the relative error. We recall that the PSNR, measured in dB, is 
defined as follows:
\[
\text{PSNR} = 
10 \log_{10} \Biggl( \dfrac{\text{MPV}^2}{\mathbf{e}_{\text{MSE}}} \Biggr)
\]
where MPV is the maximum pixel value and $\mathbf{e}_{\text{MSE}}$ is the mean squared error.
We stress that both PSNR and relative error are evaluated inside the ROI 
$S$ only. This is consistent with the motivation of ROI CT, that aims to 
recover the image inside the ROI only. The corresponding reconstructed 
images are reported in Figures \ref{fig:recIm} and \ref{fig:recEx}.

We investigated the performance of the algorithm for the regularization parameter $\lambda =5 \cdot 10^{\ell}$, 
$\ell={-4}, {-3}, \ldots, 1$;  
for the TV parameter $\rho$ we sampled the values $10^{-2}$, $10^{-1}$ and 
$1$. We considered both the explicit and the implicit formulations of the 
objective function by designing different versions of the SGP method and 
by exploiting the early stopping technique. 
As far as the Tikhonov-like regularization term concerns, we also considered 
a one-level undecimated Daubechies 4 wavelet transform in place of 
classical shearlets. However, due to room constraint, only shearlet-based 
results are reported.  In all experiments we found 
the latter to outperform the former on a visual basis, even if the figures of merit are comparable.

\begin{figure}[t]
\centering
\makebox[0pt][c]{%
\begin{tabular}{@{}cccc@{}}
$\gamma = 0.5 N$ & $\gamma = 0.3 N$ & $\gamma = 0.15 N$  \\
\includegraphics[scale=0.28]{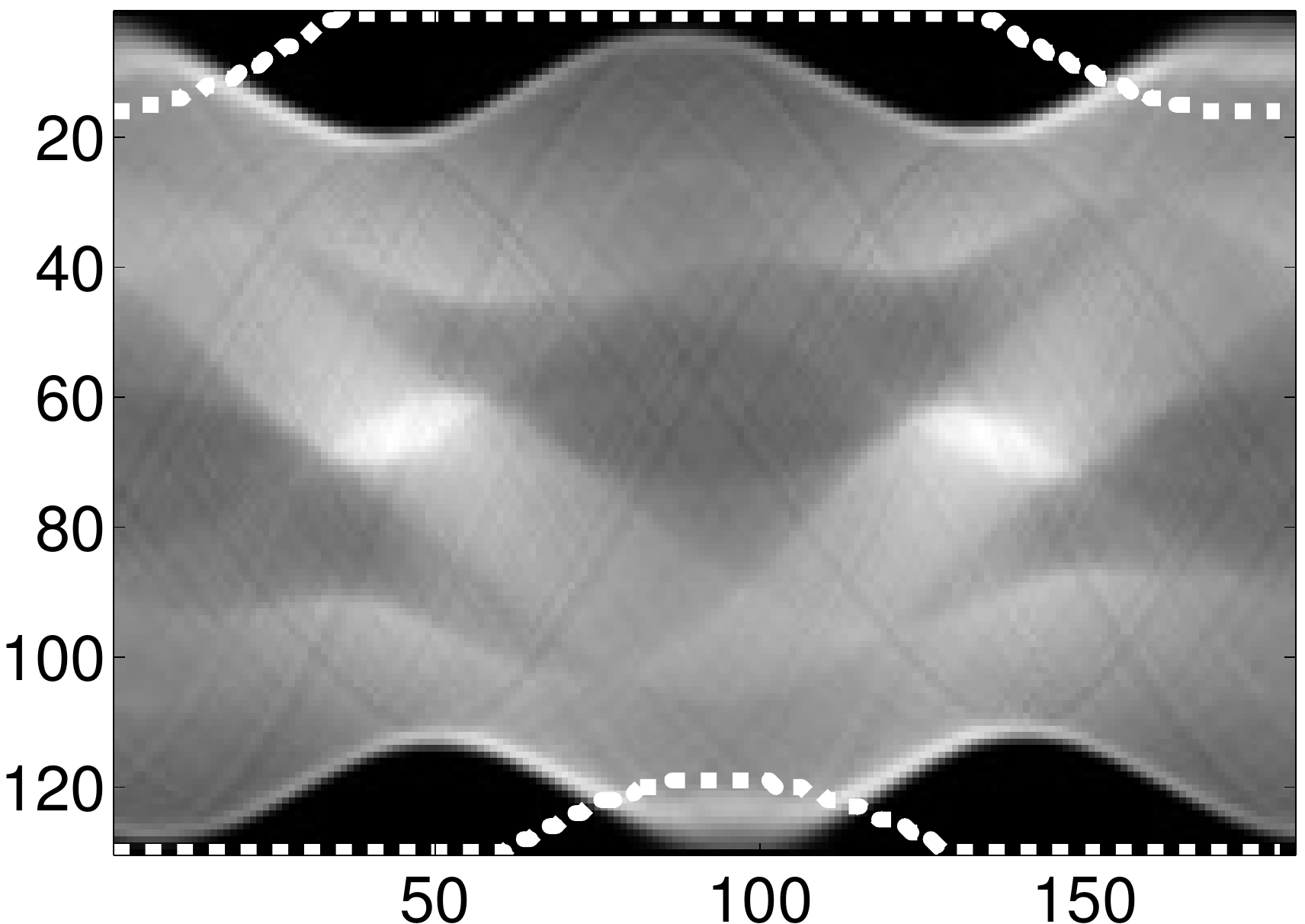}
& \includegraphics[scale=0.28]{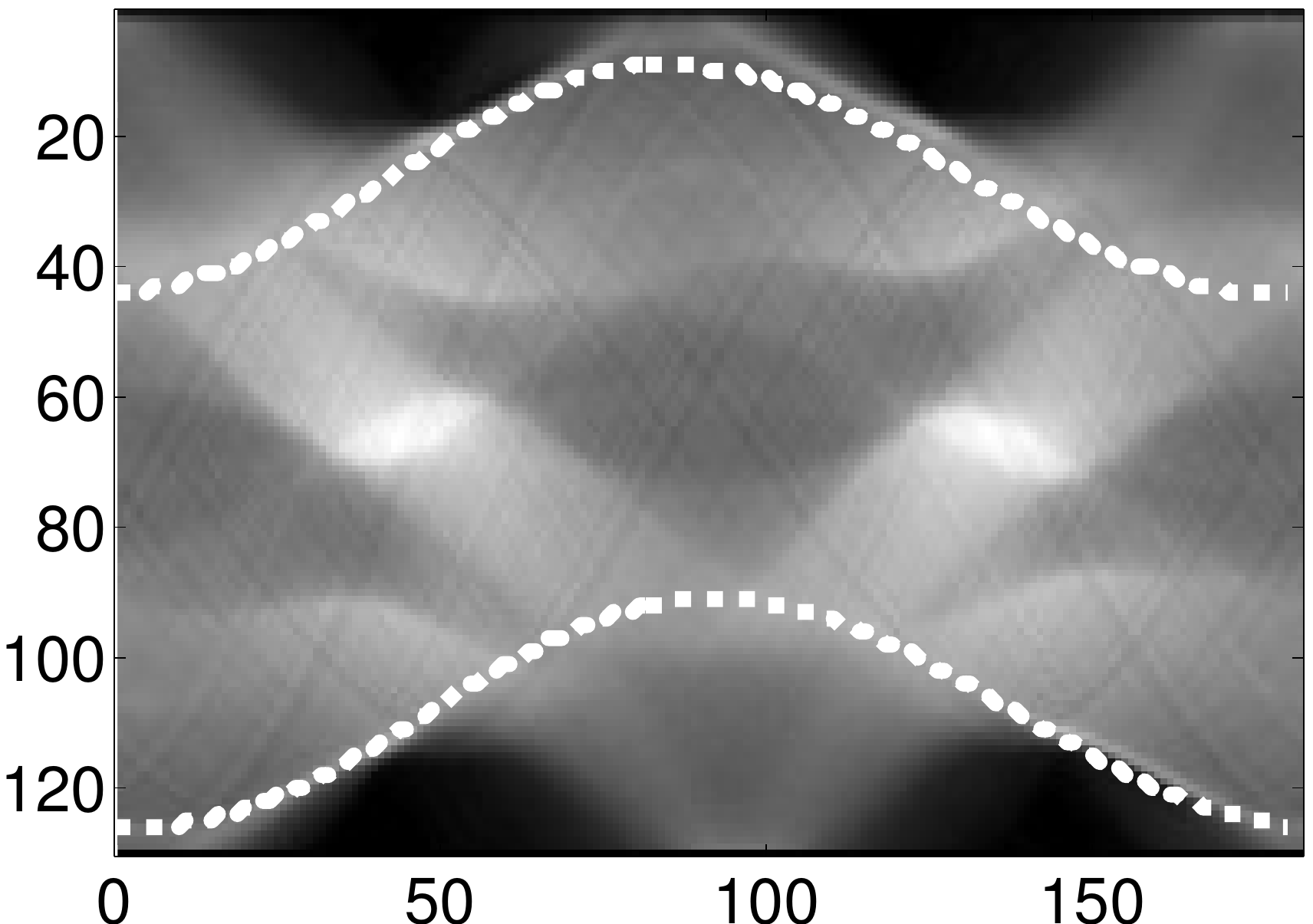}
& \includegraphics[scale=0.315]{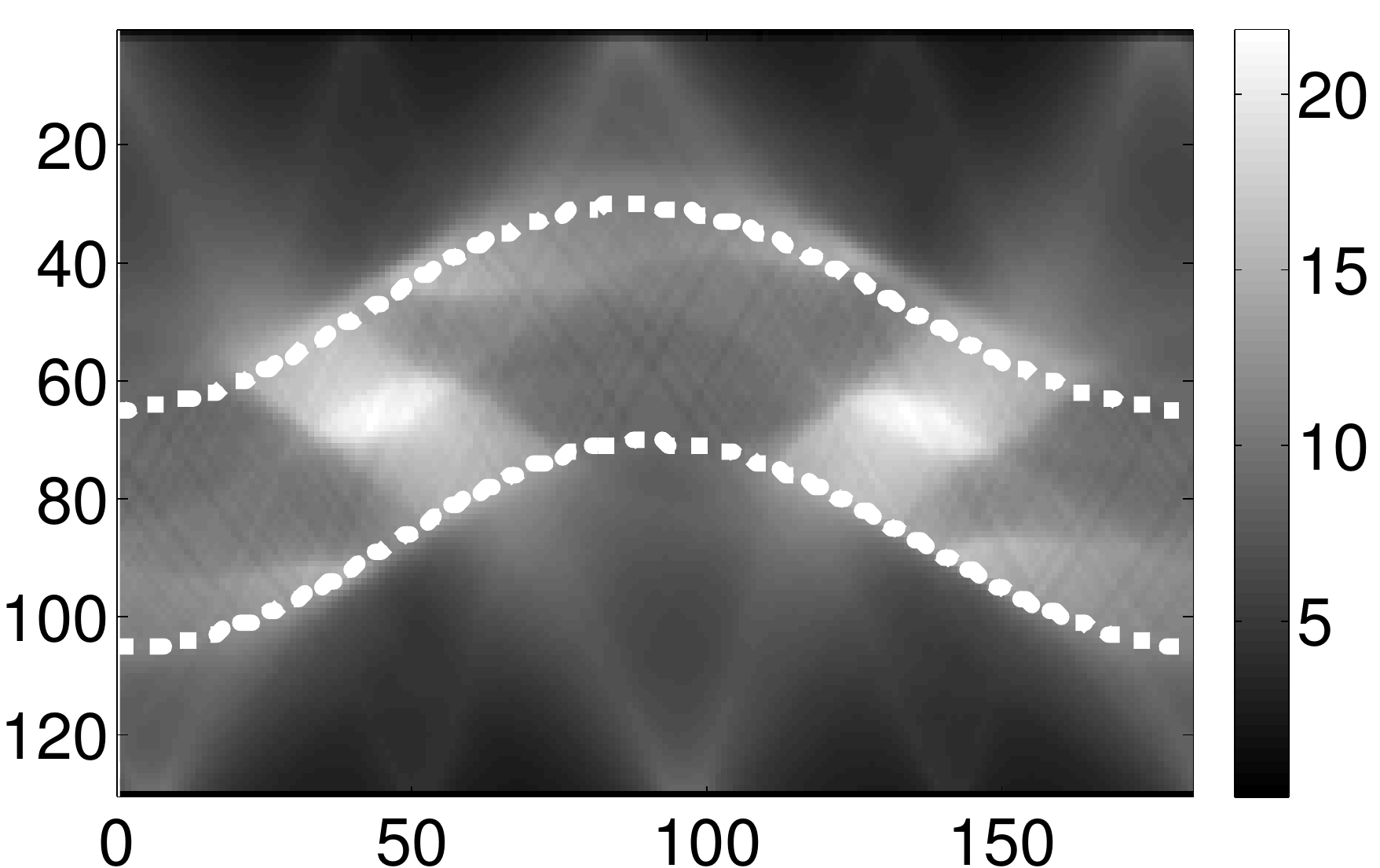} \\
(a) & (b) & (c) \\
\includegraphics[scale=0.28]{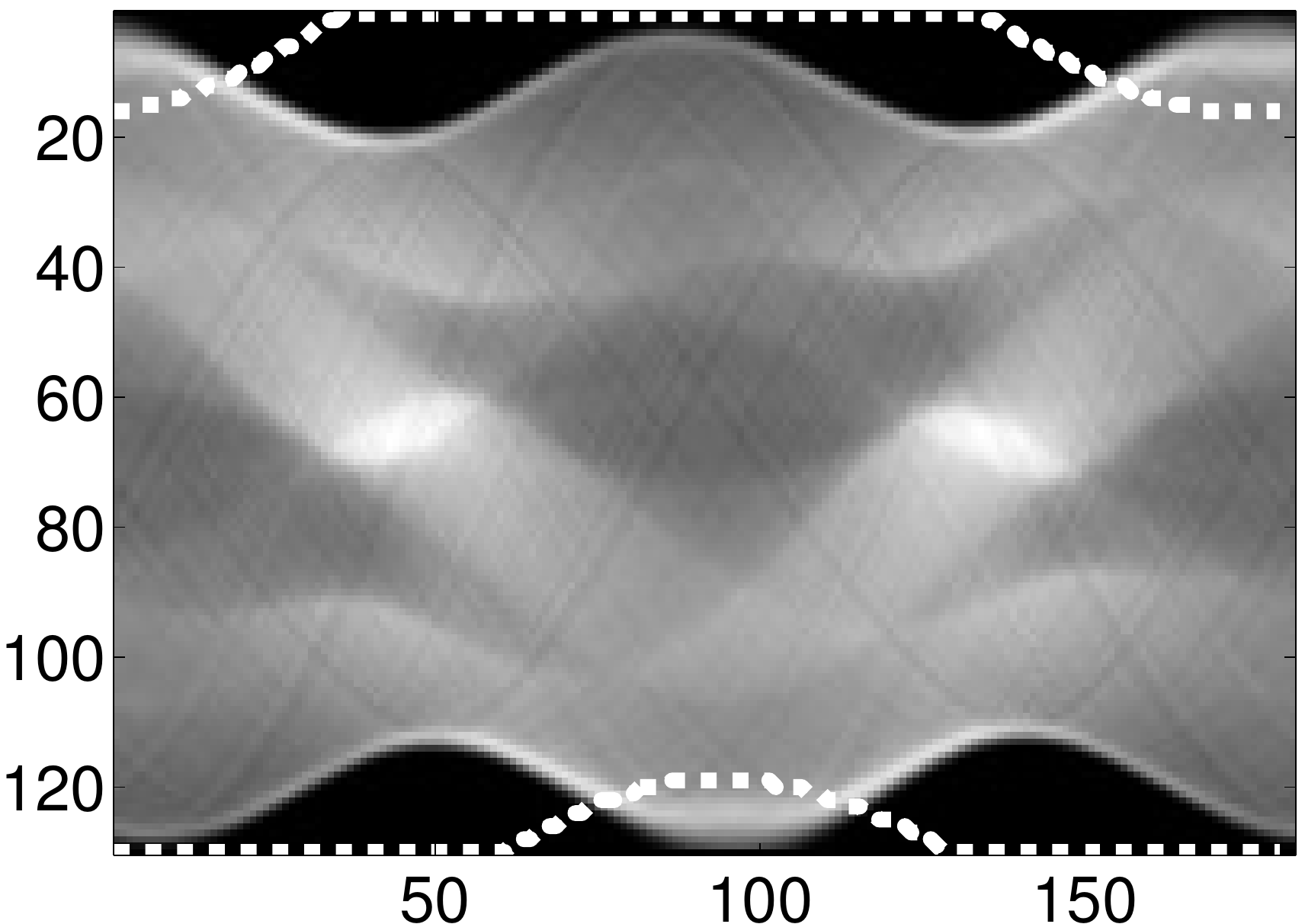}
& \includegraphics[scale=0.28]{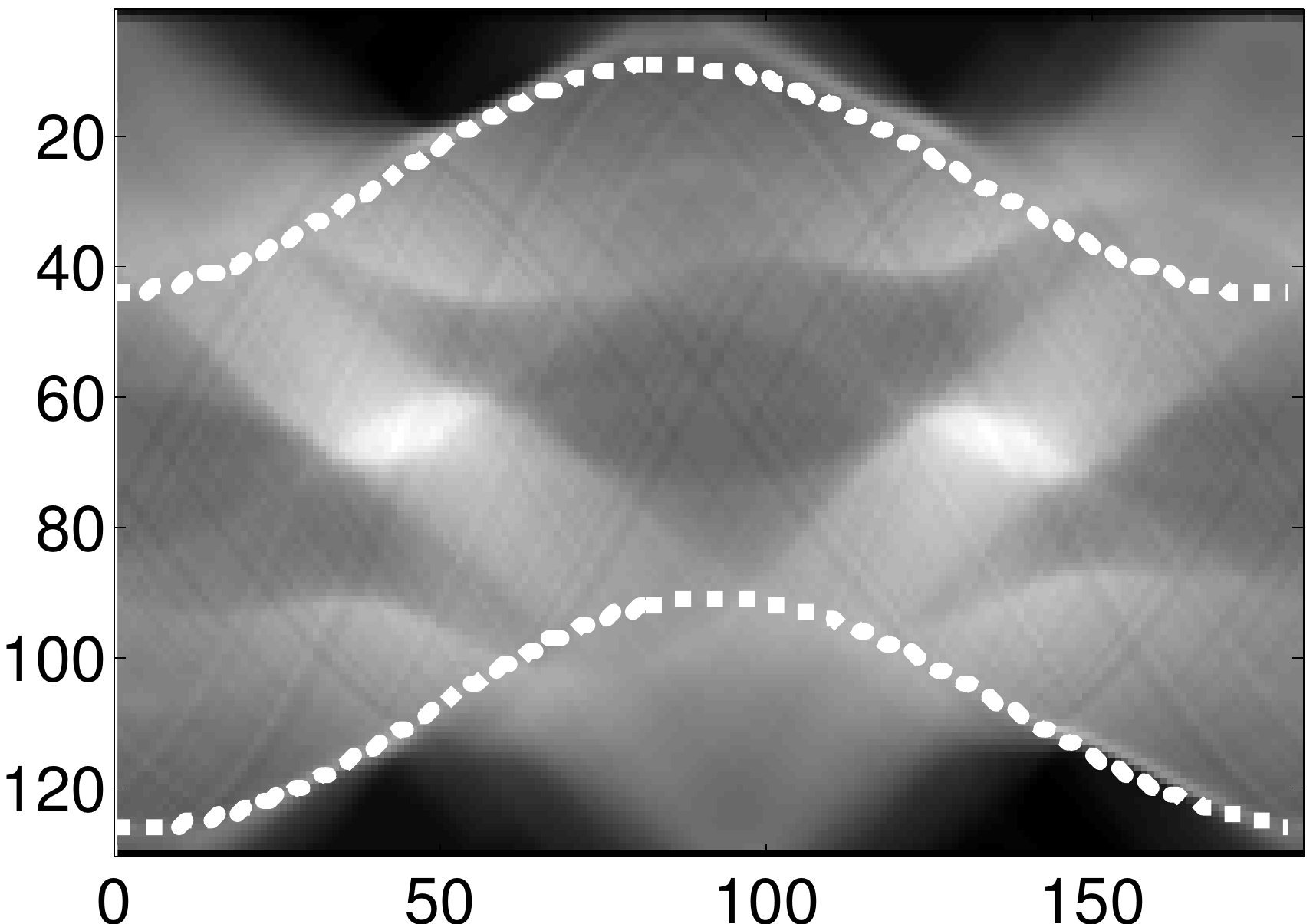}
& \includegraphics[scale=0.315]{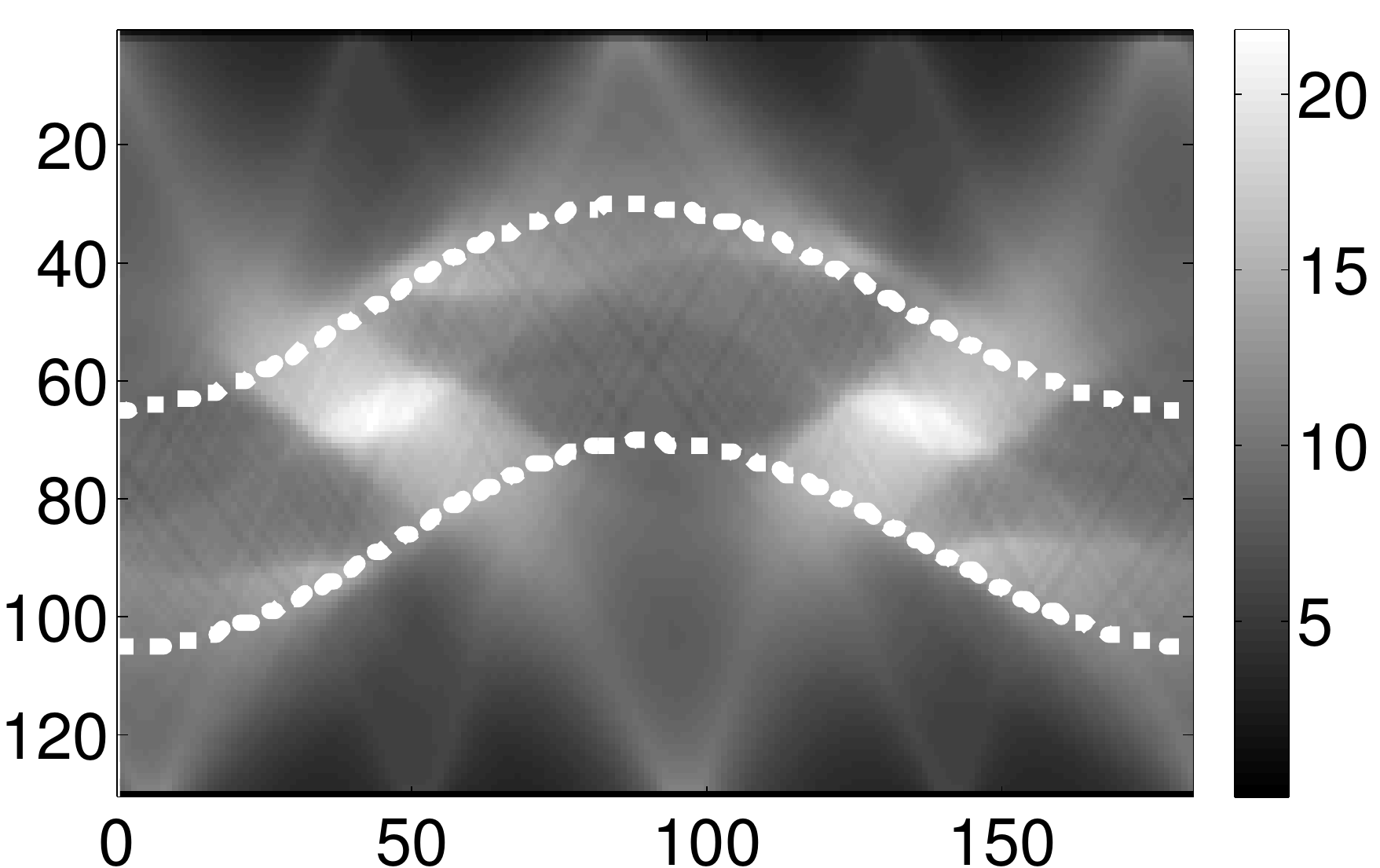} \\
(d) & (e) & (f) 
\end{tabular}}
\caption{Optimal reconstructions of the sinogram of the Shepp-Logan phantom for explicit formulation with decreasing 
\textit{radii}: $\gamma=0.5 N$ for (a) and (d), $\gamma=0.3 N$ for (b) and (e), $\gamma=0.15 N$ for (c) and (f). 
First row: shearlets and TV. Second row: just TV.
}
\label{fig:recSin} 
\end{figure}

Figures \ref{fig:recIm} and \ref{fig:recEx} report the reconstructed 
images obtained for the implicit and the explicit formulations, respectively, 
and for decreasing \emph{radii}. In details, in both figures, 
every row contains the 
reconstruction obtained with a different approach. Referring to Figure 
\ref{fig:recIm}, in the first and in the 
second row, the objective function reads as in \eqref{eq:ImplicitTV} and 
\eqref{eq:Implicit}, respectively. In the third row, the objective function is 
obtained by dropping the 
Tikhonov-like regularization term of \eqref{eq:ImplicitTV}. In the fourth row, 
both the regularization and the TV terms are dropped and the reconstruction 
is obtained by exploiting the early stopping technique. 
The same holds true for Figure \ref{fig:recEx}, where the referring equation 
for the first two rows are \eqref{eq:ExplicitTV} and \eqref{eq:Explicit}, 
respectively.
Finally, in both figures, every column
contains the reconstructions obtained for different \textit{radii}: 
$\gamma=0.5 N$ in the first column, $\gamma=0.3 N$ in the second one 
and $\gamma=0.15 N$ in the final one. 
In all of the pictures, the ROI is identified with a dashed white circle.
Figure \ref{fig:recSin} collects the reconstructed sinograms obtained for 
the explicit formulation. We only report the reconstructions corresponding to 
the first and the third row of Figure \ref{fig:recEx}, respectively. The circular ROI in the image domain corresponds to two 
sinusoidal dashed white lines in the projection domain. Finally, we report that 
further experiments with radius $\gamma$ equal to $0.25 N$, $0.2 N$ and 
$0.1 N$ were performed: the corresponding figures of merit and 
reconstruction images do not appear here. 

\begin{table}[t]
\centering
\makebox[0pt][c]{%
\begin{tabular}{@{}c|r|rrl@{\ }|@{\ }rrl@{\ }|@{\ }rrl@{}}
\multicolumn{2}{c|}{}& \multicolumn{3}{c}{$\gamma = 0.5 N$} & \multicolumn{3}{c}{$\gamma = 0.3 N$} & \multicolumn{3}{c}{$\gamma = 0.15 N$}  \\
\multicolumn{1}{c}{}&& value & iter & param & value & iter & param & value & iter & param  \\
\hline
\multirow{5}{*}{\rotatebox{90}{PSNR}}%
&SGP + Sh + TV  & 45.69 & 1140 & $\lambda=5\cdot 10^{-4}$ 
& 40.35 & 3168 & $\lambda=5\cdot 10^{-4}$
& 34.63 & 604 & $\lambda=5\cdot 10^{-4}$ \\
&& & & $\rho=0.1$ & & & $\rho=1$ & & & $\rho=1$  \\
&SGP + Sh & 41.50 & 1054 & $\lambda=5\cdot 10^{-4}$ 
& 34.86 & 1495 & $\lambda=5\cdot 10^{-4}$
& 28.46 & 2189 & $\lambda=5\cdot 10^{-4}$ \\
&SGP + TV & 49.59 & 1642 & $\rho=0.01$ 
& 41.51 & 2933 & $\rho=1$
& 39.55 & 556 & $\rho=1$ \\
&SGP & 45.01 & 7000 &  & 37.00 & 550 &  & 33.29 & 7000 &  \\
\hline
\multirow{5}{*}{\rotatebox{90}{Rel.\ err.}}%
&SGP + Sh + TV  & 0.02 & 1140 & $\lambda=5\cdot 10^{-4}$ 
& 0.10 & 3168 & $\lambda=5\cdot 10^{-4}$
& 0.50 & 604 & $\lambda=5\cdot 10^{-4}$ \\
&& & & $\rho=0.1$ & & & $\rho=1$ & & & $\rho=1$  \\
&SGP + Sh & 0.04 & 1054 & $\lambda=5\cdot 10^{-4}$ 
& 0.19 & 1495 & $\lambda=5\cdot 10^{-4}$
& 1.03 & 2189 & $\lambda=5\cdot 10^{-4}$ \\
&SGP + TV & 0.02 & 1642 & $\rho=0.01$ 
& 0.09 & 2933 & $\rho=1$
& 0.29 & 556 & $\rho=1$ \\
&SGP & 0.02 & 7000 &  & 0.15 & 550 &  & 0.59 & 7000 & \\
\hline
\end{tabular}}
\caption{Optimal results for the \emph{explicit} formulation of the objective function. Corresponding reconstructed images are reported in Figure 
\ref{fig:recEx}.
}
\label{tab:explicit}
\end{table}

In Figure \ref{fig:recIm}, we observe that, when the 
\textit{radius} is as large as half the number of pixels of the image edge 
(panel (a), (d), (g) and (j)), all the reconstructions are good.
The best approach for the objective function is 
to combine the Tikhonov-like regularization with the TV term (panel (a)), 
even if the error rate has the same magnitude for all possible approaches. 
Notice that, for all possible formulations, the values for the parameters 
$\lambda$ and $\rho$ are quite small: they are always the smallest ones of 
the investigated ranges.
For smaller \textit{radii} (panels (b), (e), (h) and (k) with 
$\gamma=0.3 N$), we notice that the fundamental structures are detected 
by all the approaches, even if the TV-free approaches exhibits 
some checkerboard effects. The best reconstruction, with respect to both 
the figures of merit, is obtained by dropping only the 
Tikhonov-like regularization term (from the formulation \eqref{eq:ImplicitTV}) 
and by setting to a medium value the TV regularization parameter $\rho$. 
Notice that the pure TV approach gains almost one order of magnitude
for both figures of merit, with respect to the other approaches.
The same considerations hold for the last case too. Indeed, when the 
\textit{radius} decreases (panels (c), (f), (i) and (l) with 
$\gamma=0.15 N$), the fundamental structures are still detected quite well in 
all the reported images, but also the checkerboard effect is still visible in the 
TV-free reconstructions. Even with this smaller \textit{radius}, the optimal 
reconstruction is achieved for the same formulation of the objective function, 
and for the same value of the regularization parameter as the larger 
\textit{radius} $\gamma=0.3 N$. The relative error 
has the same magnitude as in the case $\gamma=0.3 N$, gaining almost 
one order of magnitude with respect to the other approaches.

In Figure \ref{fig:recEx}, we observe that the fundamental structures are consistently well detected in all reconstructions. The checkerboard 
effect persists in the TV-free reconstructions (second and fourth row), even if 
it is barely visible for the largest \textit{radius}. The corresponding 
sinogram reconstructions reported in Figure \ref{fig:recSin} are accurate 
inside the ROI and therein they show no artifacts. This 
is clearly consistent with the good reconstruction of the corresponding 
images. The best approach for all \textit{radii}, with respect to the figures of merit, is 
obtained by exploiting the purely TV-based approach. However, the largest 
radius (panel (g)) needs a small value for the regularization parameter 
($\rho=0.01$) while the smaller ones (panels (h)-(i)) ask for a more 
severe TV regularization with $\rho=1$.
While for $\gamma=0.5 N$ the error rate (the PSNR, resp.) has the same 
magnitude for all possible approaches, for smaller \emph{radii} the purely 
TV-based approach gains one order of magnitude.
Notice that, even if figure (f) exhibits a definitely large relative error 
($100\%$), on a visual basis the reconstruction is acceptable, or even good. 
Such a large relative error is probably due to the checkerboard effect (sub- or under-estimation of the pixel intensity). 

In conclusion, observe that the pure TV approach, both in the implicit and in 
the explicit cases, is the one that requires more iterations to converge. 
As far as the early stopping approach concerns, when the number of iterations equals 7000, no semi-convergence property is 
observed, so the maximum number of iterations has been reached.

\section{Conclusions} \label{sec:conclusions}
In this paper, we presented a numerical assessment of the solution of the ROI CT 
problem via a first-order iterative minimization method. We considered
two different objective functions, each with a variable level of 
regularization aiming at ensuring stable reconstruction from truncated data.
Numerical experiments have shown that our reconstruction algorithms are
very satisfactory for all versions of the objective function; the main structures 
are recovered very accurately, with minimal ring artifacts. 
Using noiseless projection data, the purely TV-based approach performs slightly better than the shearlets-plus-TV approach, at least for smaller ROI
\textit{radii} (\ie{}, when the ROI is fully inside the field of view). We 
conjecture that this behavior might be dependent on the phantom features 
(which is piecewise constant) and may not  hold for more general data. 
We expect that, using more realistic sinograms and including noise,
the contribution of the shearlet term will become more relevant for the 
regularization. In fact, preliminary tests by the authors (not reported here) 
show that when projections are corrupted by noise the shearlet-based regularized 
reconstruction outperforms the one based solely on the TV term.
Another important feature of our approach is that the performance of our 
ROI reconstructions algorithm is not sensitive to the location and size of the 
ROI and performs consistently well also for rather small
ROI sizes, using both formulations of the objective function.
In future work, we plan to consider more realistic phantoms, examine a different 
formulation of the optimization problem as indicated in section 
\ref{sec:contsetup}, and extend this method to the 3D case.

\section{Acknowledgments}     

T.A.B., S.B. and G.Z. are supported by the Italian national research project 
FIRB2012, grant n.~RBFR12M3AC, and by the local research project 
FAR2014 of the University of Ferrara.
T.A.B. is supported by the Young Researchers Fellowship 2014 of the 
University of Ferrara. Also INdAM-GNCS is acknowledged. 
D.L. acknowledges partial support of NSF-DMS 1320910. 
%
%
\bibliography{BibDemIMACS,BibRoiCtIMACS}

\begin{thebibliography}{32}
\expandafter\ifx\csname natexlab\endcsname\relax\def\natexlab#1{#1}\fi
\providecommand{\bibinfo}[2]{#2}
\ifx\xfnm\relax \def\xfnm[#1]{\unskip,\space#1}\fi
\bibitem[{Barzilai and Borwein(1988)}]{Barzilai88}
\bibinfo{author}{J.~Barzilai}, \bibinfo{author}{J.M. Borwein},
  \bibinfo{title}{Two point step size gradient methods},
  \bibinfo{journal}{{IMA} J. Numer. Anal.} \bibinfo{volume}{8}
  (\bibinfo{year}{1988}) \bibinfo{pages}{141--148}.
\bibitem[{Beister et~al.(2012)Beister, Kolditz and Kalender}]{Beister12}
\bibinfo{author}{M.~Beister}, \bibinfo{author}{D.~Kolditz},
  \bibinfo{author}{W.~Kalender}, \bibinfo{title}{Iterative reconstruction
  methods in x-ray ct}, \bibinfo{journal}{Phys. Med.} \bibinfo{volume}{28}
  (\bibinfo{year}{2012}) \bibinfo{pages}{94--108}.
\bibitem[{Bonettini et~al.(2009)Bonettini, Zanella and Zanni}]{Bonettini09}
\bibinfo{author}{S.~Bonettini}, \bibinfo{author}{R.~Zanella},
  \bibinfo{author}{L.~Zanni}, \bibinfo{title}{A scaled gradient projection
  method for constrained image deblurring}, \bibinfo{journal}{Inverse Problems}
  \bibinfo{volume}{25} (\bibinfo{year}{2009}) \bibinfo{pages}{015002}.
\bibitem[{Bubba et~al.(2015)Bubba, Labate, Zanghirati, Bonettini and
  Goossens}]{Bubba15}
\bibinfo{author}{T.~Bubba}, \bibinfo{author}{D.~Labate},
  \bibinfo{author}{G.~Zanghirati}, \bibinfo{author}{S.~Bonettini},
  \bibinfo{author}{B.~Goossens}, \bibinfo{title}{Shearlet-based regularized
  {ROI} reconstruction in fan beam computed tomography}, in:
  \bibinfo{booktitle}{{SPIE} Optics \& Photonics, Wavelets And Applications
  XVI}, volume \bibinfo{volume}{9597}, \bibinfo{address}{San Diego, CA, USA},
  p. \bibinfo{pages}{95970K}.
\bibitem[{Clackdoyle and Defrise(2010)}]{clackdoyle:roi}
\bibinfo{author}{R.~Clackdoyle}, \bibinfo{author}{M.~Defrise},
  \bibinfo{title}{Tomographic reconstruction in the 21st century.
  region-of-interest reconstruction from incomplete data},
  \bibinfo{journal}{IEEE Signal Processing} \bibinfo{volume}{60}
  (\bibinfo{year}{2010}) \bibinfo{pages}{60--80}.
\bibitem[{Clackdoyle et~al.(2004)Clackdoyle, Noo, Guo and
  Roberts}]{Clackdoyle04}
\bibinfo{author}{R.~Clackdoyle}, \bibinfo{author}{F.~Noo},
  \bibinfo{author}{J.~Guo}, \bibinfo{author}{J.~Roberts},
  \bibinfo{title}{Quantitative reconstruction from truncated projections in
  classical tomography}, \bibinfo{journal}{IEEE Trans. Nuclear Science}
  \bibinfo{volume}{51} (\bibinfo{year}{2004}) \bibinfo{pages}{2570--2578}.
\bibitem[{Dai and Fletcher(2005)}]{Dai05}
\bibinfo{author}{Y.H. Dai}, \bibinfo{author}{R.~Fletcher}, \bibinfo{title}{On
  the asymptotic behaviour of some new gradient methods},
  \bibinfo{journal}{Math. Programming} \bibinfo{volume}{103}
  (\bibinfo{year}{2005}) \bibinfo{pages}{541--559}.
\bibitem[{{De Man} and Basu(2004)}]{DeMan04}
\bibinfo{author}{B.~{De Man}}, \bibinfo{author}{S.~Basu},
  \bibinfo{title}{Distance-driven projection and backprojection in three
  dimensions}, \bibinfo{journal}{Physics in Medicine and Biology}
  \bibinfo{volume}{7} (\bibinfo{year}{2004}) \bibinfo{pages}{2463--2475}.
\bibitem[{Easley et~al.(2008)Easley, Labate and Lim}]{ELL08}
\bibinfo{author}{G.R. Easley}, \bibinfo{author}{D.~Labate},
  \bibinfo{author}{W.~Lim}, \bibinfo{title}{Sparse directional image
  representations using the discrete shearlet transform},
  \bibinfo{journal}{Appl. Comput. Harmon. Anal.} \bibinfo{volume}{25}
  (\bibinfo{year}{2008}) \bibinfo{pages}{25--46}.
\bibitem[{Fletcher(2001)}]{Fletcher01}
\bibinfo{author}{R.~Fletcher}, \bibinfo{title}{On the Barzilai-Borwein method},
  \bibinfo{type}{Technical Report} \bibinfo{number}{NA/207}, Technical Report,
  Department of Mathematics, University of Dundee, \bibinfo{address}{Dundee,
  UK}, \bibinfo{year}{2001}.
\bibitem[{Frassoldati et~al.(2008)Frassoldati, Zanghirati and
  Zanni}]{Frassoldati08}
\bibinfo{author}{G.~Frassoldati}, \bibinfo{author}{G.~Zanghirati},
  \bibinfo{author}{L.~Zanni}, \bibinfo{title}{New adaptive stepsize selections
  in gradient methods}, \bibinfo{journal}{J. Industrial and Management Optim.}
  \bibinfo{volume}{4} (\bibinfo{year}{2008}) \bibinfo{pages}{299--312}.
\bibitem[{Goossens et~al.(2014)Goossens, Labate and Bodmann}]{Goossens14}
\bibinfo{author}{B.~Goossens}, \bibinfo{author}{D.~Labate},
  \bibinfo{author}{B.~Bodmann}, \bibinfo{title}{Region-of-interest computed
  tomography by regularity-inducing convex optimization},
  \bibinfo{journal}{submitted}  (\bibinfo{year}{2014}).
\bibitem[{Grippo et~al.(1986)Grippo, Lampariello and Lucidi}]{Grippo86}
\bibinfo{author}{L.~Grippo}, \bibinfo{author}{F.~Lampariello},
  \bibinfo{author}{S.~Lucidi}, \bibinfo{title}{A nonmonotone line-search
  technique for newtonÕs method}, \bibinfo{journal}{{SIAM} J. Numer. Anal.}
  \bibinfo{volume}{23} (\bibinfo{year}{1986}) \bibinfo{pages}{707--716}.
\bibitem[{Guo and Labate(2007)}]{GuoLab:2007}
\bibinfo{author}{K.~Guo}, \bibinfo{author}{D.~Labate},
  \bibinfo{title}{Optimally sparse multidimensional representation using
  shearlets}, \bibinfo{journal}{SIAM J. Math. Analysis} \bibinfo{volume}{39}
  (\bibinfo{year}{2007}) \bibinfo{pages}{298--318}.
\bibitem[{Guo et~al.(2009)Guo, Labate and Lim}]{Guo09}
\bibinfo{author}{K.~Guo}, \bibinfo{author}{D.~Labate}, \bibinfo{author}{W.Q.
  Lim}, \bibinfo{title}{Edge analysis and identification using the continuous
  shearlet transform}, \bibinfo{journal}{Applied and Computational Harmonic
  Analysis} \bibinfo{volume}{27} (\bibinfo{year}{2009})
  \bibinfo{pages}{24--46}.
\bibitem[{Guo et~al.(2006)Guo, Labate, Lim, Weiss and Wilson}]{Guo06}
\bibinfo{author}{K.~Guo}, \bibinfo{author}{D.~Labate}, \bibinfo{author}{W.Q.
  Lim}, \bibinfo{author}{G.~Weiss}, \bibinfo{author}{E.~Wilson},
  \bibinfo{title}{Wavelets with composite dilations and their {MRA}
  properties}, \bibinfo{journal}{Applied and Computational Harmonic Analysis}
  \bibinfo{volume}{20} (\bibinfo{year}{2006}) \bibinfo{pages}{202--236}.
\bibitem[{Hamelin et~al.(2010)Hamelin, Goussard, Dussault, Cloutier, Beaudoin
  and Soulez}]{Hamelin_2010}
\bibinfo{author}{B.~Hamelin}, \bibinfo{author}{Y.~Goussard},
  \bibinfo{author}{J.~Dussault}, \bibinfo{author}{G.~Cloutier},
  \bibinfo{author}{G.~Beaudoin}, \bibinfo{author}{G.~Soulez},
  \bibinfo{title}{Design of iterative roi transmission tomography
  reconstruction procedures and image quality analysis},
  \bibinfo{journal}{Medical Physics} \bibinfo{volume}{37}
  (\bibinfo{year}{2010}) \bibinfo{pages}{4577--4589}.
\bibitem[{Herman and Lent(1976)}]{Herman76}
\bibinfo{author}{G.T. Herman}, \bibinfo{author}{A.~Lent},
  \bibinfo{title}{Iterative reconstruction algorithms},
  \bibinfo{journal}{Computers in Biology and Medicine} \bibinfo{volume}{6}
  (\bibinfo{year}{1976}) \bibinfo{pages}{273--294}.
\bibitem[{Hestenes and Stiefel(1952)}]{Hestenes52}
\bibinfo{author}{M.~Hestenes}, \bibinfo{author}{E.~Stiefel},
  \bibinfo{title}{Methods of conjugate gradients for solving linear systems},
  \bibinfo{journal}{Journal of Research of the National Bureau of Standards}
  \bibinfo{volume}{6} (\bibinfo{year}{1952}) \bibinfo{pages}{409--436}.
\bibitem[{Kutyniok and Labate(2012)}]{KL12_book}
\bibinfo{author}{G.~Kutyniok}, \bibinfo{author}{D.~Labate},
  \bibinfo{title}{Shearlets: Multiscale Analysis for Multivariate Data},
  \bibinfo{publisher}{Springer}, \bibinfo{year}{2012}.
\bibitem[{Natterer(2001)}]{natterer:tomography}
\bibinfo{author}{F.~Natterer}, \bibinfo{title}{The Mathematics of Computerized
  Tomography}, \bibinfo{publisher}{SIAM: Society for Industrial and Applied
  Mathematics}, \bibinfo{year}{2001}.
\bibitem[{Natterer and Wubbeling(2001)}]{natterer:imagerec}
\bibinfo{author}{F.~Natterer}, \bibinfo{author}{F.~Wubbeling},
  \bibinfo{title}{Mathematical Methods in Image Reconstruction},
  \bibinfo{publisher}{SIAM: Society for Industrial and Applied Mathematics},
  \bibinfo{year}{2001}.
\bibitem[{Noo et~al.(2004)Noo, Clackdoyle and Pack}]{noo04}
\bibinfo{author}{F.~Noo}, \bibinfo{author}{R.~Clackdoyle},
  \bibinfo{author}{J.~Pack}, \bibinfo{title}{A two-step hilbert transform
  method for 2d image reconstruction}, \bibinfo{journal}{Physics in Medicine
  and Biology} \bibinfo{volume}{49} (\bibinfo{year}{2004})
  \bibinfo{pages}{3903--3923}.
\bibitem[{Noo et~al.(2002)Noo, Defrise, Clackdoyle and Kudo}]{noo:imrec}
\bibinfo{author}{F.~Noo}, \bibinfo{author}{M.~Defrise},
  \bibinfo{author}{R.~Clackdoyle}, \bibinfo{author}{H.~Kudo},
  \bibinfo{title}{Image reconstruction from fan-beam projections on less than a
  short scan}, \bibinfo{journal}{Physics in Medicine and Biology}
  \bibinfo{volume}{47} (\bibinfo{year}{2002}) \bibinfo{pages}{2525--2546}.
\bibitem[{R.Clackdoyle and Noo(2004)}]{clackdoyle:alcif}
\bibinfo{author}{R.Clackdoyle}, \bibinfo{author}{F.~Noo}, \bibinfo{title}{A
  large class of inversion formulae for the 2-d radon transform of functions of
  compact support}, \bibinfo{journal}{Inverse Problems} \bibinfo{volume}{20}
  (\bibinfo{year}{2004}) \bibinfo{pages}{1281--1291}.
\bibitem[{Shepp and Vardi(1982)}]{Shepp82}
\bibinfo{author}{L.~Shepp}, \bibinfo{author}{Y.~Vardi}, \bibinfo{title}{Maximum
  likelihood reconstruction for emission tomography}, \bibinfo{journal}{IEEE
  Trans Med Imaging} \bibinfo{volume}{1} (\bibinfo{year}{1982})
  \bibinfo{pages}{113--122}.
\bibitem[{Vogel(2002)}]{Vogel02}
\bibinfo{author}{C.~Vogel}, \bibinfo{title}{Computational methods for inverse
  problems}, \bibinfo{publisher}{{SIAM} Philadelphia}, \bibinfo{year}{2002}.
\bibitem[{Yan et~al.(2010)Yan, Tian, Zhu, Qin, Dai, Yang, Dong and
  Wu}]{Guorui:katsevichGPU}
\bibinfo{author}{G.~Yan}, \bibinfo{author}{J.~Tian}, \bibinfo{author}{S.~Zhu},
  \bibinfo{author}{C.~Qin}, \bibinfo{author}{Y.~Dai},
  \bibinfo{author}{F.~Yang}, \bibinfo{author}{D.~Dong},
  \bibinfo{author}{P.~Wu}, \bibinfo{title}{Fast {K}atsevich algorithm based on
  {GPU} for helical cone-beam computed tomography},
  \bibinfo{journal}{Information Technology in Biomedicine, IEEE Transactions
  on} \bibinfo{volume}{14} (\bibinfo{year}{2010}) \bibinfo{pages}{1053--1061}.
\bibitem[{Yang et~al.(2010)Yang, Yu, Jiang and Wang}]{YangYuJiang_2010}
\bibinfo{author}{J.~Yang}, \bibinfo{author}{H.~Yu}, \bibinfo{author}{M.~Jiang},
  \bibinfo{author}{G.~Wang}, \bibinfo{title}{High-order total variation
  minimization for interior tomography}, \bibinfo{journal}{Inverse Problems}
  \bibinfo{volume}{26} (\bibinfo{year}{2010}) \bibinfo{pages}{035013}.
\bibitem[{Zhang and Zeng(2007)}]{zeng:iterative}
\bibinfo{author}{B.~Zhang}, \bibinfo{author}{G.~Zeng}, \bibinfo{title}{Two
  dimensional iterative region of iterest reconstruction from truncated
  projection data}, \bibinfo{journal}{Medical Physics} \bibinfo{volume}{34}
  (\bibinfo{year}{2007}) \bibinfo{pages}{935--944}.
\bibitem[{Ziegler et~al.(2008)Ziegler, Nielsen and Grass}]{Ziegler2008}
\bibinfo{author}{A.~Ziegler}, \bibinfo{author}{T.~Nielsen},
  \bibinfo{author}{M.~Grass}, \bibinfo{title}{{Iterative reconstruction of a
  region of interest for transmission tomography}}, \bibinfo{journal}{Medical
  Physics} \bibinfo{volume}{35} (\bibinfo{year}{2008})
  \bibinfo{pages}{1317--1327}.
\bibitem[{Zou et~al.(2005)Zou, Pan and Sidky}]{zou:image}
\bibinfo{author}{Y.~Zou}, \bibinfo{author}{X.~Pan}, \bibinfo{author}{E.~Sidky},
  \bibinfo{title}{Image reconstruction in regions-of-interest from truncated
  projections in a reduced fan-beam scan}, \bibinfo{journal}{Phys. Med. Biol.}
  \bibinfo{volume}{50} (\bibinfo{year}{2005}) \bibinfo{pages}{13--28}.

\end{thebibliography}

\end{document}